\newcount\mgnf\newcount\tipi\newcount\tipoformule\newcount\greco

\tipi=2          
\tipoformule=0   


\global\newcount\numsec
\global\newcount\numfor
\global\newcount\numtheo
\global\advance\numtheo by 1

\def\senondefinito#1{\expandafter\ifx\csname#1\endcsname\relax}

\def\SIA #1,#2,#3 {\senondefinito{#1#2}%
\expandafter\xdef\csname #1#2\endcsname{#3}\else
\write16{???? ma #1,#2 e' gia' stato definito !!!!} \fi}

\def\etichetta(#1){(\veroparagrafo.\veraformula)%
\SIA e,#1,(\veroparagrafo.\veraformula) %
\global\advance\numfor by 1%
\write15{\string\FU (#1){\equ(#1)}}%
\write16{ EQ #1 ==> \equ(#1) }}

\def\letichetta(#1){\veroparagrafo.\verotheo
\SIA e,#1,{\veroparagrafo.\verotheo}
\global\advance\numtheo by 1
\write15{\string\FU (#1){\equ(#1)}}
\write16{ Sta \equ(#1) == #1 }}

\def\tetichetta(#1){\veroparagrafo.\veraformula 
\SIA e,#1,{(\veroparagrafo.\veraformula)}
\global\advance\numfor by 1
\write15{\string\FU (#1){\equ(#1)}}
\write16{ tag #1 ==> \equ(#1)}}

\def\FU(#1)#2{\SIA fu,#1,#2 }

\def\etichettaa(#1){(A\veroparagrafo.\veraformula)%
\SIA e,#1,(A\veroparagrafo.\veraformula) %
\global\advance\numfor by 1%
\write15{\string\FU (#1){\equ(#1)}}%
\write16{ EQ #1 ==> \equ(#1) }}

\def\BOZZA{
\def\alato(##1){%
 {\rlap{\kern-\hsize\kern-1.4truecm{$\scriptstyle##1$}}}}%
\def\aolado(##1){%
 {
{
 \rlap{\kern-1.4truecm{$\scriptstyle##1$}}}}}
 }

\def\alato(#1){}
\def\aolado(#1){}

\def\veroparagrafo{\number\numsec}
\def\veraformula{\number\numfor}
\def\verotheo{\number\numtheo}

\def\Eq(#1){\eqno{\etichetta(#1)\alato(#1)}}
\def\eq(#1){\etichetta(#1)\alato(#1)}
\def\leq(#1){\leqno{\aolado(#1)\etichetta(#1)}}
\def\teq(#1){\tag{\aolado(#1)\tetichetta(#1)\alato(#1)}}
\def\Eqa(#1){\eqno{\etichettaa(#1)\alato(#1)}}
\def\eqa(#1){\etichettaa(#1)\alato(#1)}
\def\eqv(#1){\senondefinito{fu#1}$\clubsuit$#1
\write16{#1 non e' (ancora) definito}%
\else\csname fu#1\endcsname\fi}
\def\equ(#1){\senondefinito{e#1}\eqv(#1)\else\csname e#1\endcsname\fi}

\def\Lemma(#1){\aolado(#1)Lemma \letichetta(#1)}%
\def\Theorem(#1){{\aolado(#1)Theorem \letichetta(#1)}}%
\def\Proposition(#1){\aolado(#1){Proposition \letichetta(#1)}}%
\def\Corollary(#1){{\aolado(#1)Corollary \letichetta(#1)}}%
\def\Remark(#1){{\noindent\aolado(#1){\bf Remark \letichetta(#1).}}}%
\def\Definition(#1){{\noindent\aolado(#1){\bf Definition 
\letichetta(#1)$\!\!$\hskip-1.6truemm}}}
\def\Example(#1){\aolado(#1) Example \letichetta(#1)$\!\!$\hskip-1.6truemm}

\def\include#1{
\openin13=#1.aux \ifeof13 \relax \else
\input #1.aux \closein13 \fi}

\openin14=\jobname.aux \ifeof14 \relax \else
\input \jobname.aux \closein14 \fi
\openout15=\jobname.aux

\let\EQ=\Eq


{\count255=\time\divide\count255 by 60 \xdef\hourmin{\number\count255}
        \multiply\count255 by-60\advance\count255 by\time
   \xdef\hourmin{\hourmin:\ifnum\count255<10 0\fi\the\count255}}

\def\oramin{\hourmin }

\def\data{\number\day/\ifcase\month\or january \or february \or march \or april
\or may \or june \or july \or august \or september
\or october \or november \or december \fi/\number\year;\ \oramin}

\newcount\pgn \pgn=1
\def\foglio{\number\numsec:\number\pgn
\global\advance\pgn by 1}
\def\foglioa{A\number\numsec:\number\pgn
\global\advance\pgn by 1}

\footline={\rlap{\hbox{\copy200}}\hss\tenrm\folio\hss}

\def\TIPIO{
\font\setterm=amr7 
\def \settepunti{\def\rm{\fam0\setterm}
\textfont0=\setterm   
\normalbaselineskip=9pt\normalbaselines\rm }\let\nota=\settepunti}

\def\TIPITOT{
\font\twelverm=cmr12
\font\twelvei=cmmi12
\font\twelvesy=cmsy10 scaled\magstep1
\font\twelveex=cmex10 scaled\magstep1
\font\twelveit=cmti12
\font\twelvett=cmtt12
\font\twelvebf=cmbx12
\font\twelvesl=cmsl12
\font\ninerm=cmr9
\font\ninesy=cmsy9
\font\eightrm=cmr8
\font\eighti=cmmi8
\font\eightsy=cmsy8
\font\eightbf=cmbx8
\font\eighttt=cmtt8
\font\eightsl=cmsl8
\font\eightit=cmti8
\font\sixrm=cmr6
\font\sixbf=cmbx6
\font\sixi=cmmi6
\font\sixsy=cmsy6
\font\twelvetruecmr=cmr10 scaled\magstep1
\font\twelvetruecmsy=cmsy10 scaled\magstep1
\font\tentruecmr=cmr10
\font\tentruecmsy=cmsy10
\font\eighttruecmr=cmr8
\font\eighttruecmsy=cmsy8
\font\seventruecmr=cmr7
\font\seventruecmsy=cmsy7
\font\sixtruecmr=cmr6
\font\sixtruecmsy=cmsy6
\font\fivetruecmr=cmr5
\font\fivetruecmsy=cmsy5
\textfont\truecmr=\tentruecmr
\scriptfont\truecmr=\seventruecmr
\scriptscriptfont\truecmr=\fivetruecmr
\textfont\truecmsy=\tentruecmsy
\scriptfont\truecmsy=\seventruecmsy
\scriptscriptfont\truecmr=\fivetruecmr
\scriptscriptfont\truecmsy=\fivetruecmsy
\def \eightpoint{\def\rm{\fam0\eightrm}
\textfont0=\eightrm \scriptfont0=\sixrm \scriptscriptfont0=\fiverm
\textfont1=\eighti \scriptfont1=\sixi   \scriptscriptfont1=\fivei
\textfont2=\eightsy \scriptfont2=\sixsy   \scriptscriptfont2=\fivesy
\textfont3=\tenex \scriptfont3=\tenex   \scriptscriptfont3=\tenex
\textfont\itfam=\eightit  \def\it{\fam\itfam\eightit}%
\textfont\slfam=\eightsl  \def\sl{\fam\slfam\eightsl}%
\textfont\ttfam=\eighttt  \def\tt{\fam\ttfam\eighttt}%
\textfont\bffam=\eightbf  \scriptfont\bffam=\sixbf
\scriptscriptfont\bffam=\fivebf  \def\bf{\fam\bffam\eightbf}%
\tt \ttglue=.5em plus.25em minus.15em
\setbox\strutbox=\hbox{\vrule height7pt depth2pt width0pt}%
\normalbaselineskip=9pt
\let\sc=\sixrm  \let\big=\eightbig  \normalbaselines\rm
\textfont\truecmr=\eighttruecmr
\scriptfont\truecmr=\sixtruecmr
\scriptscriptfont\truecmr=\fivetruecmr
\textfont\truecmsy=\eighttruecmsy
\scriptfont\truecmsy=\sixtruecmsy }\let\nota=\eightpoint}

\newfam\msbfam   
\newfam\truecmr  
\newfam\truecmsy 
\newskip\ttglue
\ifnum\tipi=0\TIPIO \else\ifnum\tipi=1 \TIPI\else \TIPITOT\fi\fi

\def\a{\alpha}
\def\b{\beta}

\def\g{\gamma}

\def\s{\sigma}
\def\t{\tau}
\def\th{\theta}

\def\z{\zeta}
\def\o{\omega}
\def\D{\Delta}
\def\L{\Lambda}
\def\G{\Gamma}
\def\O{\Omega}
\def\S{\Sigma}

\def\E{{I\kern-.25em{E}}}
\def\N{{I\kern-.25em{N}}}
\def\M{{I\kern-.25em{M}}}
\def\R{{I\kern-.25em{R}}}
\def\Z{{Z\kern-.425em{Z}}}
\def\1{{1\kern-.25em\hbox{\rm I}}}
\def\eu{{1\kern-.25em\hbox{\sm I}}}

\def\C{{I\kern-.64em{C}}}
\def\P{{I\kern-.25em{P}}}
\def\eop{{ \vrule height7pt width7pt depth0pt}\par\bigskip}



\def\AA{{\cal A}}

\def\CC{{\cal C}}

\def\LL{{\cal L}}

\def\SS{{\cal S}}
\def\TT{{\cal T}}

\def\WW{{\cal W}}
\def\VV{{\cal V}}

\def\LL{{\cal L}}
\def\XX{{\cal X}}

\def\RR{{\cal R}}

\def\chap #1#2{\line{\ch #1\hfill}\numsec=#2\numfor=1}

\def\sqr#1#2{{\vcenter{\vbox{\hrule height.#2pt
     \hbox{\vrule width.#2pt height#1pt \kern#1pt
   \vrule width.#2pt}\hrule height.#2pt}}}}


\newcount\foot
\foot=1
\def\note#1{\footnote{${}^{\number\foot}$}{\ftn #1}\advance\foot by 1}
\def\tag #1{\eqno{\hbox{\rm(#1)}}}
\def\frac#1#2{{#1\over #2}}

\def\text#1{\quad{\hbox{#1}}\quad}

\def\proof{{\noindent\pr Proof: }}

\def\remark{\noindent{\bf Remark: }}
\def\thanks{\noindent{\bf Aknowledgements: }}
\font\pr=cmbxsl10


\font\ch=cmbx12
\font\ftn=cmr8

\font\it=cmti10
\font\bf=cmbx10
\font\sm=cmr7

%
\catcode`\X=12\catcode`\@=11
\def\n@wcount{\alloc@0\count\countdef\insc@unt}
\def\n@wwrite{\alloc@7\write\chardef\sixt@@n}
\def\n@wread{\alloc@6\read\chardef\sixt@@n}
\def\crossrefs#1{\ifx\alltgs#1\let\tr@ce=\alltgs\else\def\tr@ce{#1,}\fi
   \n@wwrite\cit@tionsout\openout\cit@tionsout=\jobname.cit 
   \write\cit@tionsout{\tr@ce}\expandafter\setfl@gs\tr@ce,}
\def\setfl@gs#1,{\def\@{#1}\ifx\@\empty\let\next=\relax
   \else\let\next=\setfl@gs\expandafter\xdef
   \csname#1tr@cetrue\endcsname{}\fi\next}
\newcount\sectno\sectno=0\newcount\subsectno\subsectno=0\def\r@s@t{\relax}
\def\resetall{\global\advance\sectno by 1\subsectno=0
  \gdef\firstpart{\number\sectno}\r@s@t}
\def\resetsub{\global\advance\subsectno by 1
   \gdef\firstpart{\number\sectno.\number\subsectno}\r@s@t}
\def\v@idline{\par}\def\firstpart{\number\sectno}
\def\l@c@l#1X{\firstpart.#1}\def\gl@b@l#1X{#1}\def\t@d@l#1X{{}}
\def\m@ketag#1#2{\expandafter\n@wcount\csname#2tagno\endcsname
     \csname#2tagno\endcsname=0\let\tail=\alltgs\xdef\alltgs{\tail#2,}%
  \ifx#1\l@c@l\let\tail=\r@s@t\xdef\r@s@t{\csname#2tagno\endcsname=0\tail}\fi
   \expandafter\gdef\csname#2cite\endcsname##1{\expandafter
     \ifx\csname#2tag##1\endcsname\relax?\else{\rm\csname#2tag##1\endcsname}\fi
    \expandafter\ifx\csname#2tr@cetrue\endcsname\relax\else
     \write\cit@tionsout{#2tag ##1 cited on page \folio.}\fi}%
   \expandafter\gdef\csname#2page\endcsname##1{\expandafter
     \ifx\csname#2page##1\endcsname\relax?\else\csname#2page##1\endcsname\fi
     \expandafter\ifx\csname#2tr@cetrue\endcsname\relax\else
     \write\cit@tionsout{#2tag ##1 cited on page \folio.}\fi}%
   \expandafter\gdef\csname#2tag\endcsname##1{\global\advance
     \csname#2tagno\endcsname by 1%
   \expandafter\ifx\csname#2check##1\endcsname\relax\else%
\fi
   \expandafter\xdef\csname#2check##1\endcsname{}%
   \expandafter\xdef\csname#2tag##1\endcsname
     {#1\number\csname#2tagno\endcsnameX}%
   \write\t@gsout{#2tag ##1 assigned number \csname#2tag##1\endcsname\space
      on page \number\count0.}%
   \csname#2tag##1\endcsname}}%
\def\m@kecs #1tag #2 assigned number #3 on page #4.%
   {\expandafter\gdef\csname#1tag#2\endcsname{#3}
   \expandafter\gdef\csname#1page#2\endcsname{#4}}
\def\re@der{\ifeof\t@gsin\let\next=\relax\else
    \read\t@gsin to\t@gline\ifx\t@gline\v@idline\else
    \expandafter\m@kecs \t@gline\fi\let \next=\re@der\fi\next}
\def\t@gs#1{\def\alltgs{}\m@ketag#1e\m@ketag#1s\m@ketag\t@d@l p
    \m@ketag\gl@b@l r \n@wread\t@gsin\openin\t@gsin=\jobname.tgs \re@der
    \closein\t@gsin\n@wwrite\t@gsout\openout\t@gsout=\jobname.tgs }
\outer\def\localtags{\t@gs\l@c@l}
\outer\def\globaltags{\t@gs\gl@b@l}
\outer\def\newlocaltag#1{\m@ketag\l@c@l{#1}}
\outer\def\newglobaltag#1{\m@ketag\gl@b@l{#1}}

\def\t@gsoff#1,{\def\@{#1}\ifx\@\empty\let\next=\relax\else\let\next=\t@gsoff
   \expandafter\gdef\csname#1cite\endcsname{\relax}
   \expandafter\gdef\csname#1page\endcsname##1{?}
   \expandafter\gdef\csname#1tag\endcsname{\relax}\fi\next}
\def\verbatimtags{\let\ift@gs=\iffalse\ifx\alltgs\relax\else
   \expandafter\t@gsoff\alltgs,\fi}
\catcode`\X=11 \catcode`\@=\active
\localtags
%
\setbox200\hbox{$\scriptscriptstyle \data $}
\global\newcount\numpunt
\hoffset=0.cm
\baselineskip=14pt  
\parindent=12pt
\lineskip=4pt\lineskiplimit=0.1pt
\parskip=0.1pt plus1pt

\hyphenation{small}

\catcode`\@=11

\centerline {\bf  Typical Gibbs configurations  for  the 1d Random Field   Ising Model  }
\centerline {\bf  with  long range interaction. 
 \footnote* 
   {\eightrm Supported by:   CNRS-INdAM  GDRE 224 GREFI-MEFI,   M.C and E.O were  supported 
by Prin07: 20078XYHYS.    
}}
\vskip.5cm
 \centerline{ 
Marzio Cassandro \footnote{$^1$}{\eightrm 
Dipartimento di Fisica,  
Universit\'a di Roma ``La Sapienza'',
P.le A. Moro, 00185 Roma, Italy. cassandro@roma1.infn.it}
\hskip.2cm 
Enza Orlandi \footnote{$^2$}{\eightrm 
Dipartimento di Matematica, Universit\'a di Roma Tre, L.go S.Murialdo 1,
00146 Roma, Italy. orlandi@mat.uniroma3.it} and 
Pierre Picco \footnote{$^3$}{\eightrm LATP, CMI, UMR 6632,  CNRS,
Universit\'e de Provence,  39 rue Frederic  Joliot Curie,  13453
Marseille Cedex 13, France. picco@cmi.univ-mrs.fr }  }

 \footnote{}{\eightrm {\eightit AMS 2000 Mathematics Subject Classification}:
Primary 60K35, secondary 82B20,82B43.}
\footnote{}{\eightrm {\eightit Key Words}: phase transition, long--range interaction, random  field. 
 }

\vskip .5cm
{ \bf Abstract}
We  study a one--dimensional Ising spin systems with ferromagnetic, long--range interaction   
decaying as $n^{-2+\a}$, $\a \in [0,\frac 12]$,  in the presence of external random fields.  
 We assume that  the random fields are  given by a collection  of  symmetric, independent, 
 identically distributed 
real  random variables, gaussian or subgaussian    with variance   $\theta$. 
We show that for temperature 	and  variance  of the randomness    small enough, 
with an overwhelming probability with respect to the random fields, the typical configurations,  
  within volumes centered at the origin whose size  grow faster than any power of $\th^{-1}$,  
are intervals of $+$ spins followed by intervals of $-$ spins whose typical length 
is $ \simeq \th^{-\frac{2}{(1-2\a)}}$ for $0\le \a<1/2$
and       $\simeq e^{\frac 1 {\th^{2}}}$ for $\a=1/2$. 

\bigskip 

\chap {1 Introduction}1
\numsec= 1 \numfor= 1 
We consider a 
one dimensional  ferromagnetic Ising model with a two body  interaction     $J(n)=n^{-2+\a}$
  where  $n$ denotes the distance of the two spins and  $ \alpha \in [0,1/2]$ tunes  the  decay of the 
interaction.  
  We  add to this  term an   external random field 
   $ h[\o]: = \{ h_i [\o], i \in \Z \} $  given 
by a collection  of independent random variables, with mean zero,  symmetrically  distributed, 
 variance  $\th$, gaussian or sub--gaussian
defined on a probability space $(\O,\S,\P)$.  
We study  the magnetization profiles that are typical for the Gibbs measure when $\th$ and the 
temperature are  suitably small;  this on a 
subspace $\Omega_1(\th)\subset \Omega$ whose probability goes to $1$ when $\th \downarrow 0$. 

A systematic and successful analysis of this model for $\th=0$ {\it i.e.} when the magnetic
 fields are absent  has been already 
accomplished more than twenty years ago 
[\rcite {Ru},\rcite{D0},\rcite{D1},\rcite{D2},\rcite{Dy},\rcite {FS},\rcite{I},\rcite {ACCN},\rcite{IN}].
In particular it has been shown that it exhibits a phase transition only for $\a\in [0,1)$.
The presence of external random fields ($\th\ne 0$) modifies this picture. In [\rcite{AW}], 
it has been proved that for $\a\in [0,1/2]$
there exits an unique infinite volume Gibbs measure {\it i.e.} there is no phase transition.
More recently in [\rcite{COP3}] it has been proved that when $\a\in (1/2, \frac{\log 3}{\log 2}-1)$ 
the situation is analogous to the three
 dimensional   short range random field Ising model [{\rcite{BK}] : for temperature and variance of 
the randomness small enough, 
 there exist at least two distinct infinite volume Gibbs states, namely the $\mu^+$ and the $\mu^-$ 
Gibbs states. The proof is based on
 the notion of contours introduced in [\rcite{FS}] but using the geometrical description 
implemented in [\rcite{CFMP}] better suited to describe the contribution of the random fields. 
A Peierls argument is obtained by using a lower bound of the deterministic part of the cost to
 erase a contour
 and controlling the stochastic part.
 
The method used in [\rcite{AW}] to prove the uniqueness of the Gibbs measure is very powerful
 and general but does not provide any insight
about the most relevant spin configurations of this measure.

In this paper we show that for temperature and variance of the randomness small enough the typical 
  configurations are intervals of $+$ spins
followed by intervals of $-$ spins whose typical length  is $\th^{-\frac{2}{(1-2\a)}} $ for 
$0\le \a<1/2$ and becomes exponentially large
in terms of $\th^{-2}$ for $\a=1/2$.
When $\th>0$ the Gibbs measures are random valued measures. We need therefore to localize 
the region in which we inspect 
the system. All our results are given uniformly for an increasing sequence of intervals, 
centered in one point, with a  diameter going to 
infinity when $\th \downarrow 0$.

The modifications induced by the presence of random fields has been already studied for
 one dimensional Kac model with range
$\g^{-1}$ [\rcite{COP1},\rcite{COPV},\rcite{OP}]. In this case for $\th$ and $\g$ 
sufficiently small the typical size is $\g^{-2}$.
The results are consistent if one recalls that the random field one dimensional
 Kac model exhibits a phase transition for $\g\downarrow 0$ and 
$\th$ sufficiently small.
In the present paper the typical size is obtained estimating suitable upper and lower bounds.
The derivation of the upper bound is similar to the one used for the Kac model [\rcite{COP1}]. 
The lower bound follows from the observation that 
small intervals can be controlled with an estimate similar to those used in [\rcite{COP3}].

 \vskip 0.5cm
 \noindent{\bf Acknowledgements}
We  are indebted to  Errico Presutti for stimulating  comments  and criticism and  
 Anton Bovier for  interesting discussions.

\medskip 
\vskip 1truecm
\chap{2 Model, notations  and main results }2
\numsec= 2
\numfor= 1
\numtheo=1 

\medskip
\noindent{\bf 2.1. The model}
\medskip
Let $(\O,\AA,\P)$ be a probability space on which we  define 
$h \equiv \{h_i\}_{i\in \Z}$, a family of
independent, identically distributed symmetric random variables.
We  assume that
each  $h_i$  is Bernoulli  distributed with
$ \P[h_i=+1]=\P[h_i=-1]=1/2$.  With minor modifications that will be mentioned 
we could also  consider  the cases  of a Gaussian random variables 
with variance $1$ or a subgaussian {\it i.e.} such that $\E[\exp(t h_0)]\le 
\exp(t^2/2)\, \forall t\in \R$.  
This property is   satisfied for example 
for  $h_0=X/a$ with  $X$  an uniform random variable on  $[-a,+a]$, $a\in \R^+$
and up to an appropriate constant by any bounded symmetric random variable,
see [\rcite{KAH}] for basic properties of sub--gaussian random variables.

The spin configurations space is $\SS\equiv  \{-1,+1\}^\Z$.
If $\s \in \SS$ and $i\in \Z$,
$\s_i$ represents the value of the spin at site $i$. 
The pair interaction among spins is given by  $J(|i-j|)$ defined by
 $$J(n)= \left \{ \eqalign {& J(1) >>1; \cr &
 \frac 1 {n^{2- \alpha} } \quad \hbox {if}\quad  n >1, 
\quad \hbox {with}\,\, \alpha \in (-\infty, 1).} \right. \Eq (AP1) $$ 
 For $\L \subseteq \Z$ we set $\SS_\L=\{-1,+1\}^\L$; its elements
are   denoted by $\s_\L$; also,
if $\s \in \SS$, $\s_\L$ denotes its restriction to $\L$. Given
$\L\subset \Z$  finite,  
 define  
 $$ H_0 (\s_\L) =  \frac 12 \sum_{(i,j) \in \L \times \L}
 J(|i-j|)(1- \s_i \s_j).    \Eq(2.1b)$$
 For   $\omega \in \O$ set  $$ G(\sigma_\Lambda)[\o]:=-\theta \sum_{i\in \L} h_i[\o] \s_i . $$
   We consider the Hamiltonian   
 given by
  the  random variable on $(\O,\AA,\P)$
$$
  H (\s_\L)[\o]=  \frac 12 \sum_{(i,j) \in \L \times \L}
 J(|i-j|)(1- \s_i \s_j) +G(\sigma_\Lambda)[\o].  
\Eq(2.1)
$$
To take into account the interaction between the spins in $\L$ and
those outside $\L$ we set for $\eta \in \SS$
$$
  W (\s_\L,\eta_{\L^c}) =  \sum_{i\in \L} \sum_{j\in \L^c} J (|i-j|) (1-
\s_i \eta_j) 
\Eq(2.2)
$$
 and denote 
$$   H^\eta (\s_\L)[\o]:=    H (\s_\L)[\o]+    W (\s_\L,\eta_{\L^c}) .
\Eq(2.1a)
$$
   In the following we drop out 
the $\o$ from the notation. 
The corresponding {\sl Gibbs measure} on the finite volume $\L$,
at inverse temperature $\b>0$ with 
boundary condition  $\eta$  is a random variable with values
on the space of probability measures on $\SS_\L$  denoted   by 
$\mu^\eta_{ \L}$  
$$
\mu^\eta_{ \L}(\s_\L)
= \frac 1{Z^\eta_{ \L}} \exp\{-\b H^\eta(\s_\L)\} \quad \quad
\s_\L \in \SS_\L,
\Eq (2.3)
$$
where  $Z^\eta_{\L}$ is the normalization factor. 
 When the  configuration  $\eta$ is taken so that $ \eta_i=\tau $, $\tau =\pm 1$  
for all $i \in \Z$ 
we denote the corresponding Gibbs measure by $ \mu^+_{ \L}$ when $\tau=1$ and 
 $ \mu^-_{ \L}$  when $\tau=-1$. By FKG inequality the infinite volume limit $\L\uparrow \Z$
of $\mu_\L^+$ and $\mu_\L^-$ exists, say $\mu^+,\mu^-$.
By a result of Aizenman and Wehr, see [\rcite {AW}], 
\footnote{*}{\eightrm   A simplified  proof of this result which avoids 
the introduction of metastates,  by applying the FKG inequalities, is given  
by Bovier,  see [\rcite {Bo}], chapter 7.   Notice  that  although  we assume 
that  the distribution of the random field  has isolated point masses, 
the result [\rcite {AW}]   still holds. },   
  when $ \a \in [0,\frac 12]$
for $\P$--almost all $\o$, $\mu^+=\mu^-$
and therefore there is an unique infinite volume Gibbs measure that will be denoted by 
 $ \mu $.

 \medskip
\noindent{\bf 2.2. Main result}
\medskip

Any spin configuration  $\s\in \{-1,+1\}^\Z$ can be described in term of runs of $+1$, i.e.  
  sequences of consecutive sites $i_1,i_1+1,i_1+2\dots\in \Z$ where $\s_k=+1, \forall k\in\{i_1,\dots\}$,
 followed by  
runs of $-1$. A run  could have length $1$. To enumerate the runs we do as follows. 
Start from the site $i=0$. Let $\s_0=\t$,  $\t\in \{-1,+1\}$  call $\LL^\t_1=\LL^\t_1(\s)$ 
the run   containing  the origin, $\LL^{-\t}_2$ the run on the right of $\LL^\t_1$ and 
$\LL^{-\t}_0$ the run on the left of $\LL^{\t}_1$.  In this way to each configuration $\s$, 
we assign in a one to one way a sign $\t=\s_0$ and a family of runs $(\LL_j^{(-1)^{j+1}\t}, i \in \Z)$.
To shorten notation we drop the $(-1)^{j+1}\t$ and write simply $(\LL_j,j\in \Z)$.

Given a volume $V\subset \Z$ and a configuration $ \s_V$,  
let $e_V= e_V (\s_V)= \sup(j\in \Z : \LL_j\subset V)$ be the index of the rightmost run contained in $V$ 
    and
$ b_V= b_V (\s_V)= \inf(j\in \Z : \LL_j\subset V)$ the index of the leftmost  run contained in $V$.
We consider the sequences of runs $(\LL_j, b_V\le j\le e_V)$.

We give,  in a volume $V$ that we choose centered at the origin, in the regime $\b$ large 
and $\th$ small,   upper bound and lower bounds on 
 the length of the runs. 

In Theorem \eqv (main1alt) we show   that for volumes   larger than any inverse power of $\th$ 
up to subdominant terms with $ \P$--probability larger than  $1- e^{-g(\theta)}$, 
where $g(\theta)$ is
 a function slowly going to infinity as $ \theta \downarrow 0$,  the typical
 configurations   have   runs   with length of order $\th^{- \frac{2}{1-2\a}}$
 when $0\le \a<1/2$.   When $ \a= \frac 12$ we 
 show in Theorem \eqv (mhalf) that  with  overwhelming $\P$--probability 
  the typical run that contains the origin is exponentially long in $\th^{-2}$.

    \vskip0.5cm \noindent 
 {\bf \Theorem   (main1alt)} {\it     Let $ \alpha \in [0, \frac 12) $ and $ \zeta = \zeta (\a)$ as 
defined in \eqv(eo1),  there exist  $\th_0= \th_0 (\a)$, $\b_0=\b_0(\a)$ 
 and constants $c_i(\a)$, such that for all $0<\th\le \th_0$, for all $\b>\b_0$  
  $$ 
\beta  \ge  \frac  {\z}  {2^8   \theta^2},
\Eq  (mainbeta)$$ 
if  $0<\a<1/2$, setting  $g(\th)=(\log \frac{1}{\th})(\log\log \frac{1}{\th})$,  
with $\P$-probability larger than $1-e^{-g(\th)}$ and with a Gibbs measure 
larger than $1-e^{-g(\th)}$
the spin configurations 
are made of runs $(\LL_j, b_V\le j\le e_V)$  satisfying 
$$
c_1(\a)\left(\log \frac{1}{\th}\right)^{-\frac{2}{1-2\a}}\left(\log\log \frac{1}{\th}\right)^{-\frac{1}{1-2\a}}
\le 
\th^{\frac{2}{1-2\a}} \big|\LL_j\big|
\le c_2(\a) (\log \frac{1}{\th})(\log\log\frac{1}{\th}),
\EQ(alt1alpha)
$$
for all $j\in\{b_V,\dots e_V\}$ where $V$ is  a volume centered at 
the origin  having  diameter
$$
{\rm diam}(V) =c_0(\a)e^{g(\th)} \left(\frac{1}{\th}\right)^{\frac{2}{1-2\a}}.
\EQ(volumetheta)
$$

If $\a=0$, $g(\th)$ has to be replaced by  
$\hat g(\th) =\log\left(\frac{\log \frac{1}{\th}}{\th}\right)$ and \eqv(alt1alpha) becomes
$$
c_1(0) 
\le \th^2 \big|\LL_i\big|\le c_2(0)\left(\log \frac{1}{\th}\right)^3
\Eq(altr1alphazero)
$$
for all $j\in \{b_{\hat V}, \dots,e_{\hat V}\}$ where $\hat V$ satisfies 
$$
{\rm diam}(\hat V)=c_0(0)e^{\hat g(\th)}\left(\frac{1}{\th}\right)^2.
\Eq(hatvu)
$$
}

The proof of Theorem \eqv(main1alt) follows from  Propositions 
\eqv(1c) and \eqv(Gi2) 
and easy estimates. 

    \vskip0.5cm \noindent 
 {\bf \Theorem   (mhalf)} 
{\it    For  $ \a=1/2 $,  there exists $\th_0$ and $\b_0$ 
 and constants $c_i$, such that for all $0<\th\le \th_0$, for all $\b>\b_0$ such that 
\eqv(mainbeta) is satisfied,   
the run that contains the origin, satisfies the inequalities
$$
\exp{\frac{c_1 }{\th^2}} \le 
|\LL_1|\le 
\exp {\frac{c_2}{\th^2}}
\EQ(upphalf)
$$
with 
$\P$-probability larger than $1-e^{-\frac{c_0}{\th^2}}$
and with a Gibbs measure larger than $1-e^{-\frac{c_0}{\th^2}}$.
}
\smallskip
\Remark(meno) {\it The results for $\a=1/2$ are less general  because  the probability 
estimates for the lower bound for $\LL_i$ are not enough to extend results on  
exponential scales. 
However the  estimates for the upper bound are true on a much larger
scale, and we have results for a lot more than one run, see \eqv(ellmax) and \eqv(veealpha).
}

\vskip 1truecm
\chap{3 The upper bound }3
\numsec= 3
\numfor= 1
\numtheo=1

 \hskip 0,5cm { Let  $I \subset \Z$  be an interval,    $\tau= \pm 1$,    denote 
 $$ 
R^\tau (I)= \{ \s \in \SS: \s_i=\tau, \forall  i \in I \} 
\Eq(de1)
$$
the set of spin configurations   equal to $\tau$ in the interval $I$  
  and 
$$ 
 R  (I):=  R^+ (I) \cup  R^- (I).
\Eq(de2)
$$
  } 

Let $ L_{\max}$ be a positive integer and $V \subset \Z$ be an interval centered at 
the origin with $|V|> L_{\max}$. Denote 
 $$
\RR(V, L_{\max})= 
\bigcup_{{\scriptstyle I\subset V}
,\,\,{\scriptstyle |I|\ge L_{\max}}} R(I),
\Eq (rig1) 
$$ 
the set of spin configurations  having at least one run  of $+1$ or $-1$ larger than $L_{\max}$ in $V$. 
The main result of this section is the following 
\vskip0.5cm \noindent 
 {\bf \Proposition  (1c)} {\it Let $ \alpha \in [0, \frac 12] $, 
there exist  positive constants $c_\a$ and $c'_\a$ 
and $\th_0=\th_0(\a)$ such that for all $\b>0$, for all
decreasing real valued function $g_1(\th)\ge 1 $ defined on $\R$ that satisfies 
$\lim_{\th\downarrow 0}g_1(\th)=\infty$ there exist  an $\O_3(\a)\subset \O$ with 
$$ 
\P[\O_3(\a)]\ge 
\cases{
1-2e^{-g_1(\th)}, &if $0\le \a<\frac{1}{2}$;\cr
1-e^{-\frac 12 e^{g_1(\th)}}, & if $\a=\frac 12 $,\cr}
\EQ(stimaprob)
$$
$$
L_{max}(\a)=\cases{ c'_\a g_1(\th) \,\,\left(\frac{1}{\th^2}\right)^{\frac{1}{1-2\a}},&if $0<\a<1/2$;\cr
c'_0 g_1(\th) \,\,\,\,\frac{1}{\th^{2}} \big(\log \frac{1}{\th} \big)^2, &if $\a=0$;\cr 
c'_{1/2} e^{g_1(\th)}\,\,
e^{\frac{3}{2} \frac{8^2}{\th^2}}
(1+\frac{8}{\th})^{3} 
,&if $\a=1/2$, \cr}
\Eq(ellmax)
$$
and  an interval   $V(\a) \subset \Z $   centered at the origin  
$$
|V(\a)|=\cases{
c'_\a e^{g_1(\th)}\,\,\left(\frac{1}{\th^2}\right)^{\frac{1}{1-2\a}},& if $0<\a<1/2$;\cr
c'_0 e^{g_1(\th)} \,\,\,\frac{1}{\th^{2}} \big(\log \frac{1}{\th}\big)^2, & if $\a=0$;\cr
c'_{1/2} e^{\frac{1}{2}\exp({g_1(\th)})} \,\,e^{\frac{8^2}{\th^2}}\left(1+\frac{8}{\th}\right)^3
,& if $\a=1/2$,\cr}  
\EQ(veealpha)
$$
so that   on $\O_3(\a)$,  uniformly with respect to $\Lambda\subset \Z$,  
$$  
\sup_{\eta}
\mu^\eta_\L \left [\RR(V(\a), L_{\max}(\a)) \right] 
\le 
\cases{ 2 e^{g_1(\th)} e^{-\b c_\a \th^{-\frac{2\a}{1-2\a}}}, &if $0<\a<1/2$;\cr
2 e^{g_1(\th)} e^{-\b c_0 \log\big(\frac{1}{\th} \log \frac{1}{\th} \big)},& if $\a=0$;\cr
e^{\frac{1}{2}\exp({g_1(\th)})}\,\,e^{-\b c_{1/2}e^{\frac{8^2}{2\th^2}}}
,& if $\a=1/2$.\cr}
\EQ(gibbsmax)
$$
}

\remark 

There are various way to choose $g_1(\th)$. 
If one is interested to get a good probability estimates in \eqv(stimaprob) and to have a 
volume $L_{max}(\a)$ not too much different from the $\th^{-\frac{2}{1-2\a}}$   
in the case $0<\a<1/2$, 
one can take for $g_1(\th)$ a slowly varying function at zero. 
Note that $g_1(\th)=(\log [1/\th] )(\log\log [1/\th])$ have some advantages :
$e^{-g_1(\th)}$ decays faster than any inverse powers of $\th^{-1}$, 
the volume $V$ grows faster than any polynomials in $\th^{-1}$  
and the 
asymptotic behavior of 
\eqv(gibbsmax) 
is unaffected.

\proof
Since  $I'\subset I$,   $R(I)\subset R(I')$ we have 
$$
\bigcup_{{\scriptstyle I\subset V}
,\,\,{\scriptstyle |I|\ge L }} R(I)
\subset
\bigcup_{{\scriptstyle I\subset V}
,\,\,{\scriptstyle |I|= L}} R(I).
\Eq(fonda007)
$$
Therefore it is enough to consider the right hand side of \eqv(fonda007) 
instead of the left hand one. 

Assume that  $I=\cup_{\ell=1}^M \D(\ell)$ where
$\D(\ell)$, $\ell\in \{1,\dots,M\}$,   are adjacent  
intervals  of length $|\D|$.   
 We denote by $\D$ a generic interval  $\D(\ell)$,    $\ell\in \{1,\dots,M\}$. 
 We start estimating $\mu_{\L}^{\eta}(R^+(\D))$. 
We   bound from below     $Z^{\eta}_\L$   by     
the sum over configurations  constrained   
   to be in $R^-(\D)$ and collect the contributions of the magnetic fields in $\D$ 
both in the numerator and in the denominator.  We obtain:
 $$
     \eqalign {  \mu^\eta_{ \Lambda}(R^+ (\D)) & \le 
   \frac  {\sum_{\s_\L} e^{ -\b H^\eta(\s_\L)[\o]} \1_{ R^+(\D)} } 
{ \sum_{\s_\L}  e^{ -\b H^\eta(\s_\L)[\o]}   \1_{R^-(\Delta)} }  \cr & 
   \le  
  e^{2 \b \th  \sum_{i\in \D} h_i[\o]}  
\sup_{ \s_{\L\setminus \D} }
 \sup_{\eta_{\L^c}} \frac {  e^{ -\b   [W (\s_\D, \s_{\L\setminus \D}) + 
 W (\s_\D, \eta_\L^c)]  }  \1_{ R^+(\D)} (\s_\D)} {   e^{ -\b[   W (\s_\D, \s_{\L\setminus \D}) 
+  W (\s_\D, \eta_\L^c)] } \1_{R^-(\Delta)}(\s_\D } \cr
& \le  e^{2 \b \th  \sum_{i\in \D} h_i[\o]}  
e^{ 2 \b    [  \sum_{i\in \D} \sum_{j\in   \D^c} J (|i-j|) ]}
\le  e^{2 \b \th  \sum_{i\in \D} h_i}  e^{ 2 \b E_\a(|\D|)}.   
 \cr
} \Eq  (E.3Ca) 
$$
where $E_\a(|\D|)$ is defined by 
$$
E_\a(|\D|)=\cases{ 2(J(1)-1)+\frac{2|\D|^\a}{\a(1-\a)}, &if $0<\a<1$;\cr
2(J(1)-1) +2\log (|\D|)+ 4 ,&if $\a=0$.\cr} 
\Eq(alfapos)
$$
Calling 
$$
\O^-_1(\D)=\big\{\o : \th\sum_{i\in \D}h_i<-2E_\a(|\D|)\big\},
\Eq(omeg1)
$$
on $\O^-_1(\D)$ we have 
$$
\sup_{\L\subset\subset \Z}\sup_{\eta} \mu^\eta_\L(R^+(\D))\le e^{-2\b E_\a(|\D|)}.
\Eq(muonomega1)
$$
Define 
$$
\O^-_2(I)=\big\{\o : \exists \ell_I^*\in \{1,\dots,M\} :
\th\sum_{i\in \D(\ell_I^*)}h_i<-2E_\a(|\D|)
\big\}. 
\EQ(omega2)
$$
On $\O^-_2(I)$ we have 
$$
R^+(I)\subset R^+(\D(\ell^*_I)), 
\Eq(triv1)
$$
therefore, by  \eqv(muonomega1),
$$
\sup_{\L\subset\subset \Z}\sup_{\eta} \mu^\eta_\L(
R^+(I))\le 
e^{-2\b E_\a(|\D|)}.
\Eq(muonomega2)
$$
Assume   $V=[-N|\D|, N|\D|]$.  We can,  then,  cover $V$ 
with   overlapping intervals $I_k=[k|\D|,M|\D|+k|\D|[$
for $k\in\{-N,\dots,(N-M)\}$. It is easy to check that for 
any interval $I$ of length $M|\D|$, $I\subset V$, 
there exists an unique $k\in \{-N,\dots,(N-M-1)\}$ such that 
$$
I\supset I_k\cap I_{k+1}.
\EQ(fonda1)
$$
Therefore  one gets 
$$
\bigcup_{I\subset V,\, |I|=M|\D|}R^+(I)\subset
\bigcup_{k=-N}^{N-M-1}
\bigcup_{ {\scriptstyle I : I_{k}\cap I_{k+1} \subset I\subset V}\atop
{\scriptstyle |I|=M|\D|}} R^+(I)
\subset \bigcup_{k=-N}^{N-M-1}
R^+(I_k\cap I_{k+1}).
\EQ(fonda2)
$$
Note that for all $k$ there are $M-1$ consecutive blocks of size $|\D|$ 
in $I_k\cap I_{k+1}$ that will be indexed by $\ell_k \in \{2,\dots, M\}$.
Define   
$$
\O^-_3(V)
=\big\{\o : \forall k\in \{-N,\dots,N-M\},  
\exists \ell^*_{k} \in \{2,\dots,M\} : 
\,\,\th \sum_{i\in \D(\ell^*_k)} h_i
<-2E_\a(|\D|)\big\}.
\Eq(omega3)
$$
If we notice that $R^+(I_k\cap I_{k+1})\subset R^+(\D(\ell^*_k))$,
it  follows from \eqv(rig1), \eqv(fonda2), and \eqv(muonomega2),
that on $\O^-_3(V)$, uniformly with respect to 
$\L\subset  \Z$ we have
$$
\sup_{\eta} \mu^\eta_\L(R^+(V,M|\D|))
\le (2N+1) e^{-2\b E_\a(|\D|)}.
\Eq (muonomega3) 
$$
Next we make  a suitable choice of the parameters
$|\D|,M,N$.   Consider first the case $0<\a<1/2$. 

\noindent Since the $h_i$ are independent symmetric random variables,  we have, see \eqv (omeg1), 
$$
\P[\O^-_1 (\Delta)]=
\frac{1}{2}\Big(1-\P\big[\big|\sum_{i\in \D} h_i\big|\le \frac{2E_\a(|\D|)}{\th}\big]\Big)
\equiv \frac{1}{2}(1-p_1),
\Eq(pomega1)
$$
  
$$
\P[\O^-_2(I)]\ge 1-\left(1-\P[\O^-_1]\right)^M=1-\left(\frac{1+p_1}{2}\right)^M,
\Eq(probomega2)
$$
see \eqv (omega2),  and,  see \eqv (omega3),
$$
\P[\O^-_3(V)] \ge 1-(2N+1)\left(\frac{1+p_1}{2}\right)^{M-1}.
\Eq(pomega3)
$$
 To estimate $p_1$, we  apply   the following estimate, see  Le Cam [\rcite{Lcam}], pg  407, which holds    
for   i.i.d. random variables,     symmetric 
 and subgaussian: 
$$
\sup_{x\in \R}\P[\sum_{i=1}^{|\D|}h_i \in [x,x+\t]]\le 
\frac{2\sqrt{\pi}}
{\sqrt{|\D|\E[1\wedge(h_1/\t)^2]}}. 
\Eq(concentration)
$$
When  $ \{ h_i, i \in \Z\} $ have    symmetric  Bernoulli distribution, 
assuming that $\t\ge 1$, one has $\E[(h_1/\t)^2\1_{|h_1|\le \t}]\ge \t^{-2}$.
For  random fields having different distribution   see   Remark \eqv(r1).

For any    $0<B<1$, take  $\D$ such that $p_1\le B<1$ and $\t=2E_a(|\D|)/\th \ge1$. 
Assuming that the second constraint holds and 
using \eqv(concentration),  to satisfy the first constraint,
it is enough that
$$
p_1\le \frac{8E_\a(|\D|)\sqrt{\pi}}
{\th \sqrt{|\D|}}\le B.
\Eq(choice1)
$$
We choose
$$
|\D|=\left(\frac{32}{B\th\a(1-\a)}\right)^{\frac{2}{1-2\a}}.
\Eq(Delta)
$$
Then it is easy to check  that there exists a $\th_0=\th_0(\a,J(1))$ but independent on $B$
such that \eqv(choice1) and $\t \ge 1$ are satisfied for all $0<\th\le \th_0$. 
Choosing 
$$
M=\frac{2g_1(\th)}{\log \frac{2}{1+B}}
\EQ(choiceM)
$$
and
$$
2N+1=e^{g_1(\th)}\frac{1+B}{2} 
\EQ(choiceN)
$$
with  $g_1(\th)$   so that    $\lim_{\th\downarrow 0}g_1(\th)=\infty$,
\eqv(stimaprob), \eqv(ellmax), \eqv(veealpha), and \eqv(gibbsmax)
are proven for $0<\a<1/2$. The actual value of $B$ affects only
the values of the constants.

\noindent 
When   $\a=0$, Le Cam inequality suggests 
$$ 
|\D|=\th^{-2}
 \left(\frac{64\sqrt{\pi}}{B} \log \th^{-1}\right)^2.
\Eq(choiceD)
$$
Taking $M$ and $N$ as in \eqv(choiceM) and \eqv(choiceN), 
one gets 
\eqv(stimaprob), \eqv(ellmax), \eqv(veealpha), and \eqv(gibbsmax).

\noindent When  $\a=1/2$
$$
\O_1(\D)=\{\o : \th \sum_{i\in \D}h_i\le -8 \sqrt{\D}\}. 
\Eq(omegahalf)
$$
Le Cam inequality is useless.  We   use the Berry-Esseen Theorem 
[\rcite{CT}] that gives 
$$
\P[\O_1(\D)]\ge \frac{1}{\sqrt{2\pi}}\int_{-\infty}^{-\frac{8}{\th}}e^{-\frac{x^2}{2}}\,dx
-\frac{C_{BE}}{\sqrt{\D}}
\EQ(BE)
$$
where  $C_{BE} \le 7.5$ is the Berry-Esseen constant.    By the lower bound  
$ \int_{-\infty}^{-y}e^{-\frac{x^2}{2}}\,dx
\ge  \frac{y}{1+y^2} e^{-\frac 12 y^2}$, 
we have
$$
\frac{1}{\sqrt{2\pi}}\int_{-\infty}^{-\frac{8}{\th}}e^{-\frac{x^2}{2}}\,dx
\ge \frac{1}{\sqrt{2\pi}}\frac{1}{1+\frac{8}{\th}} e^{-\frac{8^2}{2\th^2}}.
\Eq(normal)
$$
Choosing  $$
\D=16^2 (2\pi)\left(1+\frac{8}{\th}\right)^2e^{\frac{8^2}{\th^2}},
\EQ(deltahalf)
$$
so    that the right hand side of \eqv (BE)  is strictly positive, 
$$
M=2\sqrt{2\pi}(1+\frac{8}{\th})e^{\frac{8^2}{2\th^2}}e^{g_1(\th)},
\Eq(emmhalf)
$$
and
$$
2N+1=e^{\frac{1}{2}e^{g_1(\th)}}
\EQ(ennhalf)
$$
we get  
\eqv(stimaprob), \eqv(ellmax), \eqv(veealpha), and \eqv(gibbsmax).

\eop

\Remark(r1) {\it  
 To  apply \eqv(concentration), one needs a lower bound for  the censored variance
at $\t$ of $h_1$ which is  $\E[1\wedge(h_1/\t)^2]$. 
A simple one is $\E[(h_1/\t)^2\1_{|h_1|\le \t}]$ 
which is bounded from below by half the variance of $h_1$ times $\t^{-2}$ by taking 
$\t$ large enough. However one can also get more precise bound since the difference 
between the censored variance
and the variance can be estimated by using an exponential 
Markov inequality that can be obtained as  a consequence  of the definition of sub-gaussian. 
 When $h_i, i \in \Z $    are  normal distributed   
the bound \eqv (concentration) 
  can be easily improved to 
$$
\sup_{x\in \R}\P[\sum_{i=1}^{|\D|}h_i \in [x,x+\t]]\le \frac{\t}{\sqrt{2\pi |\Delta|}}.
\EQ(gaussconcentration)
$$
}

\medskip 
\vskip 1truecm
\chap{4  Lower bound  }4
\numsec= 4
\numfor= 1
\numtheo=1

 \vskip0.5cm

Let $ \Delta  \subset \Z$  be an interval,   $\partial \Delta = \{ i \in \Z: 
 d(i, \Delta)=1 \}$, $\tau= \pm 1$,       define  
 $$ 
\WW  ( \Delta, \tau )= \{ \s \in \SS: \s_i=\tau, \forall 
 i \in \Delta, \s_{\partial \Delta }=-\tau \}.  
\Eq (CC1)$$ 
Let  $L_{\min}$ be  a  positive integer and   $ V \subset \Z$ 
 be an interval centered at the origin, with $ |V| > L_{\min}$. We denote 
 for   $i\in V$ and $\t\in \{-1,+1\}$,
$$
\nu_i(L_{\min},\t)=\bigcup_{\D\ni i,\,|\D|\le L_{\min}} \WW(\D,\t), 
\EQ(nuitau) 
$$

$$ 
\VV (V, L_{\min})= \bigcup_{i\in V} \left[\nu_i(L_{\min},+) \cup \nu_i(L_{\min},-)\right].
\Eq (notrig1)
$$
The main result of this section is the following. 

  \vskip0.5cm \noindent  
{\bf \Proposition  (Gi2) } {\it     Let  $ \alpha \in [0, \frac 12] $,
 $\theta>0$, $ \zeta= \zeta (\a)$ as defined in \eqv (eo1). 
There exists $\th_0=\th_0(\a)$ and $ \b_0=\b_0(\a)$ such that for $0<\th<\th_0$ and 
$\b>\b_0$, for all $D>1$, for all decreasing real valued function $g_2(x)$ defined 
on $\R^+$ such that 
$\lim_{x\downarrow 0}g_2(0)=\infty$ but $\lim_{x\downarrow 0}\frac{g_2(x)}{x}=0$, if 
we denote 
$$ \bar b:= \min(\frac {\beta \zeta} 4, \frac  {\z^2}  {2^{10}  \theta^2})  
\Eq (E.8a) $$
then  there exists $\O_5(\a)\subset \O$ with
$$
\P[\O_5(\a)]\ge 
\cases{ 1-  
5\left(\bar b\right)^{\frac{2}{(1-2\a)}}e^{-(4D-1)g_2(\bar b)}
,& if $0<\a<1/2$;\cr
1-5\left(\frac{\bar b}{g_2(\bar b)}\right)^2
\left(4+\log\left[\frac{\bar b}{8g_2(\bar b)}\right]\right)^2
e^{-(4D-1) g_2(\bar b)}
,& if $\a=0$;\cr
1-e^{-g_2(\bar b)},&if $\a=1/2$.\cr
}\Eq(pomega5)
$$
For  
$$
L_{min}(\a)=
\cases{
\left(\frac{\bar b}{D g_2(\bar b)} \right)^{\frac{1}{1-2\a}}
 \left(
\frac{1}{4+\log\left(\bar b \right)^{\frac{1}{1-2\a}}}\right)^{\frac{1}{1-2\a}}
,&if $0<\a<1/2$;\cr
\frac{\bar b}{Dg_2(\bar b)} \left(4+\log\big[\frac{\bar b}{Dg_2(\bar b)}\big]\right),
&if $\a=0$;\cr
e^{\frac{\bar b}{2D}-4} ,& if $\a=1/2$,\cr
}
\EQ(ellmina)
$$
and 
$$
V_{min}(\a)=\cases{
e^{g_2(\bar b)}(\bar b)^{\frac{1}{1-2\a}},& if $0<\a<1/2$;\cr
e^{g_2(\bar b)}
\frac{\bar b}{Dg_2(\bar b)} \left(4+\log\big[\frac{\bar b}{Dg_2(\bar b)}\Big]\right),&
if $\a=0$;\cr
e^{\frac{\bar b}{2}(1-\frac{1}{D})}e^{-2g_2(\bar b)}
,&if $\a=1/2$,\cr}
\EQ(veemina)
$$
 on $\O_5(\a)$, for all    $ \L \subset \Z$ large enough,  
 $$  
\mu^+_\L  [   \VV (V_{min}(\a), L_{min}(\a))]   
\le 
\cases{
5\left(\bar b\right)^{\frac{2}{(1-2\a)}}e^{-(4D-1)g_2(\bar b)}
,&if $0<\a<1/2$;\cr
5\left(\frac{\bar b}{8D g_2(\bar b)}\right)^2
\left(4+\log\left[\frac{\bar b}{D g_2(\bar b)}\right]\right)^2
e^{-(4D-1)g_2(\bar b)}
,&if $\a=0$;\cr
e^{-g_2(\bar b)},& if $\a=1/2$.\cr}
\Eq  (E.3c)$$  

}
 \vskip0.5cm

\remark  The estimate \eqv  (E.3c)  is uniform in $\L$, therefore by  the  uniqueness of the infinite volume Gibbs measure,  
[\rcite {AW}],   Proposition  \eqv (Gi2)  holds for the infinite volume Gibbs measure $\mu$.

\proof
 Since the  boundary conditions are homogeneous equal to $+$ we    apply  the   geometrical description
 of the spin configuration  presented in  [\rcite {CFMP}].   In the following 
we will assume 
that the notions of triangles, contours, and their properties are known to the reader.
 In  Section 7   we summarize definitions and main properties  used in the proof.  
 Let    $ \TT = \{  \underline T \} $  be  the  set  of  families of   triangles  
 compatible with the chosen  $+$ boundary conditions on $\L$.   Let denote by $|T|$ the 
 {\it mass } of  the triangle $T$, {\it i.e.}  the cardinality of $ T \cap \Z$, see \eqv (mars1).  
It is convenient to identify  in  $ \underline T \in \TT$   families  of triangles
having the same mass, 
$$
\underline T =  \{\underline T^{(1)},\dots,  \underline T^{(k_{ \underline T})} \}, 
\Eq (mars10) 
$$
arranged in increasing order, 
where  $ k_{ \underline T}= \sup \{|T|\,: T \in \underline T \}  \in \N$
and for $\ell \in \{1,\dots,k_{\underline T}\}$, $\underline T^{(\ell)}$ is the family
of $n_\ell\equiv n_\ell(\underline T) \in \N$ triangles in $\underline T$ having all 
the mass $\ell$. 
By convention $n_\ell(\underline T)=0$ when there is no triangle of mass $\ell$ in $\underline T$.
We denote
$$    
|\underline T|^x = \sum_{\ell=1}^{k_{ \underline T}}  n_\ell(\underline T) \,\ell^x   
, \,\,x \in \R, \,\,x \neq 0, 
\Eq (Ma.3) 
$$    
 and
$$  
\log   |\underline T|= \sum_{\ell=1}^{k_{ \underline T}}  n_\ell(\underline T) 
(4+  \log {\ell}).   
\Eq (Ma.3z) 
$$
Let $ \Lambda \subset \Z$ be an interval large enough, $ V \subset \L$ and $L$ an integer,  $L \le |V|$. 
 Since 
$\mu^+_\L(\cup_{i\in V }\nu_i(L,-))\le \sum_{i\in V}\mu^+_\L(\nu_i(L,-))$, it is enough to
 estimate for a given $i\in V$,
$\mu^+_\L(\nu_i(L,-))$. Applying  \eqv(nuitau) one has
$$
\mu^+_\L(\nu_i(L,-))\le \sum_{\ell_0=1}^L\,\,\sum_{\D: \D\ni i, |\D|=\ell_0}
\mu^+_\L(\WW(\D,-)). 
\EQ(P1)
$$
It remains to estimate $\mu^+_\L(\WW( \Delta, -))$, for a given $i\in V$, $\Delta\ni i$ and 
$|\D|=\ell_0$.   
We denote by 
 $$ 
 \CC=\CC(\D,-) = \{ \underline T \in  \TT   \hbox {compatible  with } \WW( \Delta, -)  \}.
\EQ(calc)
 $$
A family $\underline T$ is said   compatible with the event $\WW(\Delta,-)$ if 
$\underline T$ corresponds to a spin configuration where the event $\WW(\D,-)$ occurs. 
By construction   the     families  of triangles  in $  \CC$ satisfy  
 only one of  the two    following conditions:
 
\item{$\bullet$}    there exists $ T_0 \in \CC$  so that   $ \Delta = supp (T_0)$ 

\item{ $\bullet$}       there exist  two triangles $T_{right}=T_{right}(\D)$ and $T_{left}
=T_{left}(\D)$ 
one on   the right and one on  the left of $\Delta$  that are 
  adjacent \footnote*  {\eightrm  
 We say that $T$ is adjacent to an interval $\Delta$ if   
$ 0<d(supp(T), \Delta)< 1  $. {\it i.e.} 
 $\Delta \cap  supp(T) = \emptyset$  and $T$ is the first triangle on  the right  
  or the left of $\Delta$ having the support  at distance from $\Delta$ smaller than 1.  }  
to $\Delta$.

The fact that $T_{left}$ (resp. $T_{right}$) is on the left (resp. right) of $\D$
and is adjacent to it will  be denoted by $T_{left} \triangleleft\D$,
(resp $T_{right} \triangleright \D$).
By  \eqv (Ma1)   $\ell_0={\rm dist}(T_{left},T_{right})\ge |T_{right}|\wedge |T_{left}|$, {\it i.e.}
at least  one of the two triangles  $(T_{left},T_{right})$  has  support  smaller or equal 
than $\ell_0$. 
We make the partition:
$$  
\CC =    \cup_{j=1}^3 \AA_j   
\Eq (E.1) $$
where $\AA_j=\AA_j(\D,i)$ are defined by:
 $$ \AA_1 = \{ \underline T\in  \CC : \exists T_0 \in  \underline T,    
supp (T_0)=\D \}; \Eq (do1)$$
 $$ 
\AA_2 =  \cup_{\ell=1}^{\ell_0}  \AA_2 (\ell)
\,\,{\rm with} \,\,\, \AA_2 (\ell)=   \{ \underline 
T \in   \CC  : \exists T_{left}\in  \underline T,\, T_{left}\triangleleft\, \D,
\,    |T_{left}|= \ell  \};    
\Eq (do2) $$
  $$ \AA_3 =  \cup_{\ell=1}^{\ell_0}  \AA_3 (\ell)\,\,{\rm with} \,\,\,
  \AA_3 (\ell) =  \{ \underline T \in \CC  \setminus    \AA_2 : \exists  
T_{right } \in  \underline T, \,T_{right}\triangleright \D,\,    |T_{right}|= \ell \}.  
\Eq (do3) $$ 
Any family in $\AA_1$ can be written as $(T_0, \underline T)\in \AA_1$
where $T_0\notin \underline T$.    
  We denote by  $\AA_1\setminus T_0$ the set all these $\underline T$ such that 
$(T_0,\underline T) \in \AA_1$,
    with the same meaning  we denote    
 $\AA_2(\ell)\setminus T_{left}$ and $\AA_3(\ell)\setminus T_{right}$.
We have
$$
\eqalign { 
   \mu^+_{ \Lambda}  ( \WW( \Delta, -))  =& 
\sum_{  \underline  T   \in \AA_1 \setminus T_0} \mu^+_{ \Lambda} 
 ( T_0\cup \underline   T)\1_{\{Supp(T_0)=\D\}}
\cr & + 
\sum_{\ell=1}^{\ell_0} \,\,
\sum_{T_{left}: |T_{left}|=\ell}\1_{\{T_{left} \,\triangleleft\,\, \D\}}
\sum_{  \underline  T  \in \AA_2(\ell) \setminus  T_{left} } 
\mu^+_{ \Lambda}  ( T_{left} \cup \underline   T  )  
  \cr &+\sum_{\ell=1}^{\ell_0}\,\,  
\sum_{T_{right}: |T_{right}|=\ell}\1_{\{T_{right}\, \triangleright\,\,\D \}}
\sum_{  \underline  T  \in \AA_3(\ell)
 \setminus  T_{right} } \mu^+_{ \Lambda}  ( T_{right}\cup \underline   T ) .  
} 
\Eq (N.1)$$
For any given triangle $T$, with $|T|=\ell$, recalling the definition 
of contours in   Section 7,  let 
$$ 
\AA (T)  \equiv  \AA (T, \ell)= \{  
\underline S \in  \TT: T\notin \underline S\, ;   (T, \underline S) 
\,\, \hbox {form a contour}\,;
 \,\,  \forall S \in  \underline S,\,   |S | <\ell   \}. 
\Eq (mar10) 
$$

 \vskip0.5cm \noindent    \Remark   (mar11)  All the  triangles  belonging to   $\AA (T, \ell)$  have  mass   $\ell_1 <\ell $  and    form    a contour with  $T$. 
 Notice that    triangles   $T_1$  with $    |T_1 |= \ell_1 $,    $\ell_1 <\ell $ might belong  to the same contour $\G$   of $T$ but when   we remove  the  triangles  in $\G$ different than $T$,  having support larger or equal to  $\ell$ the resulting family might not form a single contour with $T$. 
 These triangles are not in $ \AA (T, \ell)$.

  \vskip0.5cm \noindent 
  We start analyzing  the   first term on the right hand side of \eqv(N.1).
We decompose $\underline T\in \AA_1\setminus T_0$ as 
$\underline S_1\cup \underline T'$ with $\underline S_1\in \AA(T_0,\ell_0)$ and 
$\underline T' \notin \AA(T_0,\ell_0)$,  obtaining
$$
\eqalign{
\sum_{  \underline  T   \in \AA_1 \setminus T_0} \mu^+_{ \Lambda} 
 ( T_0\cup \underline   T) 
=&\sum_{\underline S_1 \sim T_0}\1_{\{\underline S_1 \in \AA(T_0,\ell_0)\}}
\sum_{\underline T'\sim (T_0\cup\underline S_1)}
\1_{\{\underline T'\notin \AA(T_0,\ell_0)\}}
\mu^+_\L(T_0\cup\underline S_1\cup\underline T')\cr
=&\sum_{\underline S_1 \sim T_0}\1_{\{\underline S_1 \in \AA(T_0,\ell_0)\}}
\mu^+_\L(T_0\cup\underline S_1)
=\sum_{\underline S_1 \in \AA(T_0,\ell_0)}
\mu^+_\L(T_0\cup\underline S_1). 
\cr
}\EQ(new1pi)
$$
Recall that $\underline S_1\sim T_0$ means that 
$\underline S\cup T_0$ is an allowed configuration of triangles.
Applying the  same  decomposition for the  remaining  two terms on the right hand side of \eqv(N.1)  we get 
$$ 
\eqalign { 
  \mu^+_{ \Lambda}  ( \WW( \Delta, -))  =&   
 \sum_{  \underline  S_1 \in \AA(T_0, \ell_0) }
 \mu^+_{ \Lambda} ( T_0\cup   \underline S_1 )  
\1_{\{supp(T_0)=\D\}}
 \cr 
& +
 \sum_{\ell=1}^{\ell_0}\,  
\sum_{T_{left}: |T_{left}|=\ell}\1_{\{T_{left}\,\,\triangleleft\,\, \D  \}}
\sum_{  \underline  S_1  \in \AA(T_{left}, \ell)  } 
\mu^+_{ \Lambda}  ( T_{left} \cup \underline  S_1 )  \cr 
&+\sum_{\ell=1}^{\ell_0}\,  \sum_{T_{right}: |T_{right}|=\ell}
\1_{\{T_{right}\,\triangleright\,\,\D  \}}
\,
\sum_{  \underline  S_1   \in \AA(T_{right}, \ell) } 
\mu^+_{ \Lambda}  ( T_{right} \cup \underline  S_1 ). 
} \Eq (LL.1)$$
We estimate separately each  term in the previous sums.  They are   all
alike $ \mu^+_{ \Lambda} (T\cup \underline S)$ with $\underline S \in \AA(T,\ell)$ see 
\eqv(mar10) and 
$|T|=\ell$. 
Recalling \eqv(mars10), we identify in  $\underline S$ the families of triangles
 having the same mass. 
By construction  we have $k_{\underline S}\in \{1,\dots,\ell-1\}$. 
We follow  an argument used in [\rcite{COP3}] which consists  of  4 steps. 
We consider first  the case $0<\a<1/2$, the case  $\a=0$ and $\a=\frac 12$ will be discussed later.

{\noindent \bf Step I}

For each  $ j= \{ 1,\dots,   k_{ \underline S}\} $ we extract 
a term  $\sum_{k=1}^{ j}   n_k(\underline S)  k^\a$ 
from  the deterministic part of the Hamiltonian, {\it i.e.} using Theorem \eqv(1CFMP),
we write  
$$  
 \eqalign { 
 \mu^+_\L (T\cup \underline S )=&  
\frac 1 { Z^+_\Lambda[\o]} \sum _{  \underline T'\sim T\cup \underline S  } 
 e^{ -\b H^+( \underline T' \cup T \cup \underline   S )[\o]}  \cr &
 \le e^{-\b\frac {\zeta} 2  (  \sum_{k=1}^{ j}   n_k(\underline S)  {k}^\a) } 
    \frac 1 { Z^+_\Lambda[\o]}  \sum_{   \underline T' \sim T\cup \underline S }  
e^{ -\b H^+_0( \underline T'\cup T \cup \underline   S  \setminus  
 (\cup_{k=1}^j \underline S^{(k)}  )  + \beta \theta  
G (\s( \underline T' \cup T \cup \underline S  ))[\o]}.  }  
\Eq (D.2aa) 
$$
We add to this list  of $k_{\underline S} $ inequalities a $k_{\underline S}+1$--th inequality   
that we get when,   after   extracting  all the terms
corresponding to $\underline S$,  we extract the term corresponding to  $T$ {\it i.e.}
$$   
 \mu^+_\L (T\cup \underline S )    
 \le e^{-\b\frac {\zeta} 2  (  \sum_{k=1}^{ k_{ \underline S} }  
 n_k(\underline S)  {k}^\a  +  \ell ^\a ) }  
    \frac 1 { Z^+_\Lambda[\o]} \sum _{  \underline T' \sim T\cup \underline S}
  e^{ -\b H^+_0(  \underline T' ) 
 + \beta \theta  G (\s(  \underline T'  \cup \underline S   \cup T ))[\o]}. 
\Eq (D.2bb) 
$$
Observing the right hand side of \eqv(D.2aa) and \eqv(D.2bb), one notes 
that the $H^+_0$ and $G$ are not evaluated at the same configuration of triangles.  
In the next step we  compensate this discrepancy by a corrective term. 

{\noindent \bf Step II}

For each $j \in  \{1,\dots, k_{\underline S} \} $      we multiply and divide  \eqv (D.2aa) by
  $$    
\sum _{  \underline T' \sim  T\cup\underline S } 
e^{ - \b H^+_0(   \underline T'\cup T\cup \underline   S   \setminus  
 (\cup_{\ell=1}^j \underline S^{(\ell)}  )  + \beta \theta  
G (\s(  \underline T'   \cup T \cup  \underline S \setminus  
 (\cup_{\ell=1}^j \underline S^{(\ell)}  ))[\o]}   
 \Eq (D.4) 
$$  
and   when $  j = k_{\underline S}+1$, see  \eqv (D.2bb)  by
$$    
\sum _{  \underline T' \sim      T\cup \underline S }
 e^{ -\b H^+_0(   \underline T')  + 
\beta \theta  G (\s(  \underline T'))[\o]}.   
\Eq (D.4a) $$ 
Setting for $ j\in \{  1,\dots,   k_{ \underline S}\} $ 
$$ 
 F_j[\o]: =   
\frac 1 \b \ln \left \{  \frac {  \sum_{\underline T'\sim T\cap \underline S}  
e^{-\beta H^+_{0} (\underline T' \cup T  \cup \underline S  \setminus   (\cup_{\ell=1}^j 
\underline S^{(\ell)}   )  + \beta \theta  G (\s(  \underline T'  \cup T   
  \cup \underline S ))[\o]} }
  {  \sum_{   \underline T'\sim    T\cup \underline S} 
  e^{-\beta H^+_{0} (\underline T' \cup T  \cup \underline S \setminus  
 (\cup_{\ell=0}^j \underline S^{(\ell)}  ) + \beta \theta 
 G  (\s (\underline T'\cup T  \cup \underline S \setminus   
(\cup_{\ell=1}^j \underline S^{(\ell)}  ))[\o]} }\right\} ,  
\Eq (L.1b) 
$$ 
 and for $j = k_{ \underline S}+1$ 
$$  
F_{ k_{ \underline S}+1} [\o] =  
 \frac 1 \b \ln \left \{  \frac 
{  \sum_{\underline T' \sim T\cup\underline S} 
 e^{-\beta H^+_{0} (\underline T' )  
+ \beta \theta  G (\s(  \underline T  \cup T  \cup \underline S ))[\o]} }
  {  \sum_{   \underline T'\sim    T\cup \underline S   } 
  e^{-\beta H^+_{0} (\underline T' ) + 
\beta \theta  G  (\s (\underline T'))[\o]} }\right\}   
\Eq (L.1ba) 
$$ 
we have the following set of inequalities: 
for   $j \in \{1, \dots, k_{ \underline T}+1 \}$ 
$$
\mu^+_\L (T \cup \underline S) 
\le
e^{-\b\frac {\zeta} 2  (  \sum_{\ell=1}^{ j}   n_\ell(\underline S)  {\ell}^\a )  
+  \b F_j[\o]}
     \mu_\L^+ ( T \cup \underline S \setminus   (\cup_{\ell=1}^j \underline S^{(\ell)} ))  
\le e^{-\b\frac {\zeta} 2  (  \sum_{k=1}^{ j}   n_k(\underline S)  {k}^\a )  +  
\b F_j[\o]} .
 \Eq (MM1) $$

{\noindent \bf Step III}

We make a partition of the probability space to take into account the 
fluctuations of the $F_i$ in \eqv(MM1). 
For   each    $ (T,\underline S) $ we write  
$$
\O=
\cup_{j=0}^{k_{\underline S}+1}  B_j , 
\Eq (P.3) 
$$  
where, recalling \eqv(Ma.3),  for $j\in \{1, \dots,  k_{ \underline S}\}$ 
$$
B_j=B_j ((T,\underline S))  
=\{ \o\, :   F_j[\o]  \le 
 \frac \zeta 4     
\sum_{k=1}^j 
n_k(\underline S) \,k^{ \alpha},   
  \hbox { \rm and for} \;
 \forall i \in \{j+1, \dots \ell_0\},
   F_i [\o]  >   \frac \zeta 4  
   \sum_{k=1}^i n_k(\underline S) 
\,k^{ \alpha}  \};  
 \Eq (P.1) 
$$    
$$
B_{ k_{ \underline S}+1}=B_{ k_{\underline S}+1} ((T,\underline S)) 
 = \left \{ \o\,:   F_{ k_{ \underline S}+1}[\o]  \le  \frac \zeta 4   
  \left (\sum_{k=1}^{k_{ \underline S}} n_k(\underline S)\, k^\a
 +  
 \ell^{ \alpha}\right ) \right \};   
\Eq (P.1b) 
$$
$$
B_0=B_0 ((T,\underline S))  =\{ \o\,:    
 \forall i \in \{1, \dots,  k_{ \underline S}+1\},  
 F_i [\o]  >   \frac \zeta 4     \sum_{k=1}^i
 n_k(\underline S)\, k^\a
\}.   
\Eq (P.1a) 
$$ 
The point  is that using exponential inequalities for Lipschitz function of
subgaussian random variables, see [\rcite{COP3}] Section 4 for details, one has : 
for all $\a \in (0,1)$
For  $0 \le  j \le  k_{\underline S} $,   
 $$ \E \left [ \1_{B_j} \right ] \le 
e^{-   \frac {\z^2 }{2^{10}\theta^2}   
 \left ( \sum_{k =j+1}^{ k_{\underline S} }  
n_k(\underline S)\, k^{2\alpha-1}  
 +  
\ell^{2\alpha-1}\right ) }.  
\Eq (E.3a)
$$
with the convention that an empty sum is zero.  For $j=k_{ \underline S}+1$ we use  
$ \E \left [ \1_{B_{k_{ \underline S}+1}} \right ] \le 1$.

{\noindent \bf Step IV}

Using \eqv(P.3), we have
$$
\E\left[\mu^+_\L (T \cup \underline S)\right] =\sum_{j=0}^{k_{\underline S}+1} 
\E\left[\mu^+_\L (T \cup \underline S)\1_{\{B_j\}}\right], 
\EQ(tr1)
$$
then,     \eqv(MM1) entails
$$
\E\left[\mu^+_\L (T \cup \underline S)\1_{\{B_j\}}\right] 
\le  e^{-\b\frac {\zeta} 2  (  \sum_{k=1}^{ j}   n_k(\underline S) \, {k}^\a )}  
\E\left[e^{\b F_j}\1_{\{B_j\}}\right] .
\Eq(triv2)
$$
Recalling \eqv(P.1),\eqv(P.1b) and \eqv(P.1a), on $B_j$ we have 
$$
F_j\le 
 \frac \zeta 4     
\sum_{k=1}^j 
n_k(\underline S)\, k^{ \alpha}
\EQ(trivb5)
$$
that gives with \eqv(triv2) and \eqv(E.3a)
$$
\E\left[\mu^+_\L (T \cup \underline S)\1_{\{B_j\}}\right] 
\le 
e^{-\b\frac {\zeta} 4  \sum_{k=1}^{ j}   n_k(\underline S) \, {k}^\a }  
e^{-   \frac {\z^2 }{2^{10}\theta^2}   
 \big( \sum_{k =j+1}^{ k_{\underline S} }  
n_k(\underline S)\, k^{2\alpha-1}  
 +   
\ell^{2\alpha-1}\big) }.  
\EQ(trivb6)
$$
Coming back to \eqv(tr1) we get 
$$
\eqalign{
\E\left[\mu^+_\L (T \cup \underline S)\right] &
    \le \sum_{j=0}^{ k_{\underline S} }   
 e^{-\frac {\beta \zeta} 4 \sum_{k=1}^j n_k(\underline S)
\,k^\alpha }  e^{-  
     \frac {\z^2} {2^{10} \theta^2  } 
\left (  \sum_{k =j+1}^{k_{\underline S}   } 
n_k(\underline S) \, k^{2\alpha-1} + \ell^{2\alpha-1}\right ) } 
+  e^{-\frac {\beta \zeta} 4  
\left ( \sum_{k=1}^{ k_{\underline S} }  
 |\underline S^{(k)}|
n_k(\underline S) \,k^\alpha    + \ell^{\alpha } \right )  } 
\cr & 
\le
(k_{\underline S} +2) e^{-\bar b
 \left (  \sum_{k =1}^{ k_{\underline S}  }
 n_k(\underline S)
\,k^{2\alpha-1} +\ell^{2\alpha-1} \right )   }, 
}
\EQ(En.10b)
$$
where 
$$
 \bar b= \min\left(\frac {\beta \zeta} 4,   \frac  {\z^2}  {2^{10}  \theta^2}\right). 
 \Eq (E.8) 
$$
 {\bf Final conclusions}
 To estimate  \eqv(P1) we take into account the  partition done in \eqv (LL.1). 
 Corresponding to  the first term in \eqv(LL.1), using \eqv(En.10b) and \eqv(do1),
 we have for each $i\in V$
$$
\eqalign{
I_1(i)&\equiv 
\sum_{\ell_0=1}^L\,\,\sum_{\D: \D\ni i, |\D|=\ell_0}
\,\,
\sum_{  \underline  S_1 \in \AA(T_0, \ell_0) }
\E\left[ \mu^+_{ \Lambda} ( T_0\cup   \underline S_1 ) \right]
\1_{\{suppT_0=\D\}}
\cr
&=
\sum_{\ell_0=1}^{L}\,\,\sum_{T_0: T_0\ni i, |T_0|=\ell_0}
\,\,\sum_{  \underline  S_1 \in \AA(T_0, \ell_0) }
\E\left[ \mu^+_{ \Lambda} ( T_0\cup   \underline S_1 ) \right]\cr
&\le 
 \sum_{\ell_0=1}^{L}\,\,\sum_{T_0: T_0\ni i, |T_0|=\ell_0}
\,\,\sum_{  \underline  S_1 \in \AA(T_0, \ell_0) }
(\ell_0 +2) e^{- \bar b
 \left (  \sum_{k =1}^{ k_{\underline S}  }
 n_k(\underline S_1)
\,k^{2\alpha-1} +\ell_0^{2\alpha-1} \right )   }.\cr
}
\Eq(PPX1)
$$
Since all  the triangles in $\underline S_1\in\AA(T_0, \ell_0)  $ 
are  smaller than $\ell_0$, we have
  $$
 \sum_{k =1}^{k_{\underline S_1}   } n_k(\underline S_1)\,k^{2\alpha-1} + \ell_0^{2\alpha-1}
\ge 
 \frac 1  {\ell_0^{1-2\a} (4+\log \ell_0)} 
 \left (  \sum_{k =1}^{k_{\underline S_1}  }n_k(\underline S_1)(4+\log k)
  + (4+\log \ell_0)\right )
\EQ(LB1)
$$
so that from \eqv (PPX1) we have 
$$
\eqalign{
I_1(i)  &\le 
 \sum_{\ell_0=1}^{L}(\ell_0 +2) \sum_{T_0: T_0\ni i, |T_0|=\ell_0}
\,\,\sum_{  \underline  S_1 \in \AA(T_0, \ell_0) }
  e^{- \bar b
 \left (   \frac 1  {\ell_0^{1-2\a} (4+\log \ell_0)} 
 \left (  \sum_{k =1}^{k_{\underline S_1}  }n_k(\underline S_1)(4+\log k)
  + (4+\log \ell_0)\right ) \right )  }
  \cr & \le
   \sum_{\ell_0=1}^{L}(\ell_0 +2) \,\sum_{\G: \G\ni i, |\G|\geq \ell_0}  e^{- \bar b
 \left (   \frac 1  {\ell_0^{1-2\a} (4+\log \ell_0)} 
 \left (  \sum_{k =1}^{k_{\underline \G}  }n_k(\underline \G)(4+\log k)
  + (4+\log \ell_0)\right ) \right )}.\cr
}
\Eq(PPc1)
$$
 Take    $D>1$ and    $ g_2(\bar b)>1$ so that 
$$
  \frac {\bar b}  {  L^{1-2\a} (4+\log L)} \ge   D g_2(\bar b).
 \Eq(hyp1)
 $$
 Applying  \eqv(E.6ab), if $ D g_2(\bar b)\ge C_0\vee 3$ 
we get 
$$
I_1(i) 
\le 
\sum_{\ell_0=1}^{L} (\ell_0+2)e^{-Dg_2(\bar b)(4+\log \ell_0)}
 \sum_{\ell_2\ge \ell_0} 2\ell_2 e^{-Dg_2(\bar b)(4+\log \ell_2)}
 \le  10  e^{-8Dg_2(\bar b)}.
\Eq(PPX3)
$$
It remains to consider the second term in \eqv(LL.1), the third term being identical. 
Using \eqv(En.10b), \eqv(do2), and \eqv(E.6ab)   
for each $i\in V$, we have
$$
\eqalign{
I_2(i)&\equiv \sum_{\ell_0=1}^L \,\,\sum_{\D : \D\ni i, \,|\D|=\ell_0}
\sum_{\ell_1=1}^{\ell_0}\,\,
\sum_{T_{left}: |T_{left}|=\ell_1}\1_{\{T_{left}\,\triangleleft\,\,\D \}}
\sum_{ \underline T\in \AA_2(\ell_1)\setminus T_{left}}
\E\left[\mu^+_\L(T_{left}\cup\underline T)\right]\cr
&\le 
 \sum_{\ell_0=1}^L\ell_0\sum_{\ell_1=1}^{\ell_0} (\ell_1+2)
e^{-Dg_2(\bar b)(4+\log \ell_1)}
\sum_{\G : \G\ni 0; |\G|\ge \ell_1}e^{-Dg_2(\bar b)(4+\log |\G|)}\cr
&\le 5 e^{-8Dg_2(\bar b)}
\sum_{\ell_0=1}^L\ell_0\le 5L^2  e^{-8Dg_2(\bar b)}.
\cr
}
\Eq(PPX4)
$$
Collecting \eqv(PPX3) and \eqv(PPX4) one gets
$$
\E\left[\mu^+_\L(\nu_i (L,-))\right]\le 20 L^2 e^{-8 Dg_2(\bar b)}.
\EQ(Re2)
$$
By  Markov inequality,  
 on a probability subset $\O_4=\O_4(L,i)$ with
 $$
 \P[\O(L,i)]\geq 1- 5Le^{-4Dg_2(\bar b )}, 
 \Eq(Omegadelta)
 $$
 one gets
 $$
 \mu^+_\L(\nu_i (L,-))\le 5L e^{-4Dg_2(\bar b )}.
 $$
 Recalling the definition of $\VV(V,\ell_0)$ see \eqv(notrig1), 
one gets that
 on a probability subset $\O_5=\O_5(V)$ with 
 $$
 \P[\O_5]\ge 1-|V|5Le^{-4D g_2(\bar b)} 
 \Eq(omegave)
 $$
 we have
 $$
 \mu^+_{\L}(\VV(V,L)\le 5|V|Le^{-4D g_2(\bar b)}. 
 \Eq(final)
 $$
 \noindent 
 {\bf Choice of the parameters}

\noindent $\bullet$  $0<\a <\frac 12.$  From  \eqv (hyp1) we take 
 $$
L\equiv  L_{\min}=
\left(\bar b \right)^{\frac{1}{1-2\a}}
 \left(4+\log\left(\bar b \right)^{\frac{1}{1-2\a}}\right)^{-\frac{1}{1-2\a}}
\left( D g_2(\bar b )\right)^{-\frac{1}{1-2\a}}.
 \EQ(ellzero)
 $$
 It is easy to check that there exists a $\th_0=\th_0(\a)$ and $\b_0$ that depend
 on $\a$ but not on $D>1$ nor  on $g_2(\bar b)\ge 1$
 such that  \eqv(hyp1) is satisfied for all $0<\th\le\th_0$ and all $\b\ge \b_0$.
 
 Then one can take the volume $V$ with a diameter similar to  \eqv(veealpha), namely
 $$
 V_{min}(\a)=e^{g_2(\bar b)}\left(\bar b\right)^{\frac{1}{1-2\a}}.
 \Eq(vmin)
 $$
 An easy computation gives \eqv(pomega5) and \eqv(E.3c).
\medskip
\noindent $\bullet$ $\a=0.$
Going back to \eqv(D.2aa), the modifications are the following :
each time a $k^\a$, respectively an $\ell^\a$, appears replace it by $(4+\log k)$, respectively 
by $(4+\log \ell)$. The event in the step III, are modified in the same way. 
The only difference comes with \eqv(E.3a) replaced by 
$$ 
\E \left [ \1_{B_j} \right ] \le e^{-   \bar b   
 \left ( \sum_{k=j+1}^{ k_{\underline S} }  n_k( \underline S)  \frac {  \left ( 4 +\log {k} 
  \right)^2  } { {k}   }  + \frac { (4+ \log  \ell_0 )^2 } {\ell_0}  \right ) }.  
\Eq (E.3azero)
$$
Then \eqv(LB1) is modified using 
$$
\frac {(4+\log k)^2}{k} \ge 
\frac {(4+\log \ell_0)}{\ell_0} (4+\log k). 
\EQ(LB01)
$$
The assumption \eqv(hyp1) becomes 
$$
\bar b \,\,\frac{4+\log L}{L}\ge  D g_2(\bar b).
\Eq(hyp1zero)
$$
Then everything but the choice of $L$ goes as before. 
Here we choose 
$$
L\equiv  L_{\min}=\frac{\bar b}{Dg_2(\bar b)} 
\left(4+\log\left[\frac{\bar b}{Dg_2(\bar b)}\right]\right)
\EQ(ell00)
$$
and it is easy to see that if $\bar b\ge D g_2(\bar b)$ then 
\eqv(hyp1zero) is satisfied.
Then as before taking
$$
V_{min}(0)=
e^{g_2(\bar b)}
\frac{\bar b}{Dg_2(\bar b)} \left(4+\log\big[\frac{\bar b}{Dg_2(\bar b)}\Big]\right)
\EQ(vmin0)
$$  
one gets \eqv(pomega5) and \eqv(E.3c) after easy estimates. 

\medskip \noindent
\noindent $\bullet$ $\a=1/2.$      \eqv(E.3a)  holds  in the following form  
$$
 \E \left [ \1_{B_j} \right ] \le e^{-   \bar b   
 \left ( 1+ \sum_{k=j+1}^{ k_{\underline S} }  n_\ell(\underline S)  
 \right ) }.  
\Eq (E.3aze)
$$
Since  $1+\sum_{k=1}^{\ell} n_k (\underline S)\ge 1$ 
the inequality \eqv(En.10b) becomes
$$
\E[\mu_\L^+(T\cup \underline S)]\le (k_{\underline S}+2)
e^{-\frac{\bar b}{2}}e^{-\frac{\bar b}{2}
 \left( 1+\sum _{k=1}^{k_{\underline S}} n_k(\underline S)\right)}.
\Eq(E.3ahalf)
$$
The  condition  \eqv(hyp1) becomes 
$$
\frac{\bar b}{2(4+\log L)}\ge D \ge C_0
\EQ(hyphalf)
$$
where $C_0$ is defined in \eqv(2CFMP).
Taking 
$$
L\equiv L_{min}=e^{\frac{\bar b}{2D}-4}
\Eq(lminhalf)
$$
one has 
$$
\E[\mu_\L^+(\nu_i(L,-))]\le 
20 
e^{+\frac{\bar b}{2D}-8} e^{-\frac{\bar b}{2}}
\le 20 e^{-\frac{\bar b}{2}(1-\frac{1}{D})}.
\EQ(Re1half)
$$
Therefore if one takes
$$
V_{\min}(1/2)=\frac{1}{20} e^{\frac{\bar b}{2}(1-\frac{1}{D}) }e^{-2g_2(\bar b)}
\EQ(vminhalf)
$$
one gets
$$
\mu_\L^+(\VV(V_{\min},L_{\min}))\le e^{-g_2(\bar b)}
\EQ(muhalf)
$$
with a $\P$--probability larger than $1-e^{-g_2(\bar b)}$.\eop

\medskip 
\vskip 1truecm
\chap{7 Appendix:  Geometrical description of the spin configurations } 7
\numsec= 7
\numfor= 1
\numtheo=1

  \vskip0.5cm

We will follow  the geometrical description of the spin configuration  presented in  [\rcite {CFMP}]
and  use the same notations.
We will consider homogeneous boundary conditions, i.e the spins in the
boundary conditions are either all $+1$ or all $-1$. Actually we will restrict ourself to
$+$ boundary conditions and  consider spin   configurations $ \s= \{ \s_i, i\in \Z \}\in \XX_+ $ 
so that $\s_i=+1$ for all $|i|$ large enough.

In one dimension an interface at $(x, x+1)$ means  $\s_{x} \s_{x+1}= -1$.   
Due to the above choice of the boundary conditions, any $\s \in \XX_+ $ has a finite, 
even number of  interfaces.   The precise location of the interface is immaterial and 
 this  fact has been  used  to choose the interface points as  follows:   For
 all $x \in \Z$ so that $(x,x+1)$ is an interface take the  location of the 
interface  to be a point  inside the interval  
$[x+\frac 12- \frac 1 {100}, x+\frac 12+ \frac 1 {100} ] $, with the property that for 
any four distinct points $r_i$,  $i=1, \dots, 4$ $|r_1-r_2| \neq |r_3-r_4|$.   This choice 
is done once for all so that the interface between $x$ and $x+1$ is uniquely fixed.     
 Draw  from each one of  these interfaces points  two 
lines  forming  respectively   an angle of $ \frac \pi 4$  and of   $ \frac 3 4 \pi  $ with 
the $\Z$ line.    We have thus a bunch of growing $\vee-$ lines each one emanating from an
 interface point.    Once two $\vee-$ lines meet, they are frozen and stop their growth. 
 The other two lines emanating  from the the same interface points are {erased}.
The $\vee-$ lines 
 emanating from others points  keep growing.  The collision of the two lines is represented
 graphically by a triangle whose basis is the line joining the two interfaces points and whose 
sides are the two segment of the $\vee-$ lines which meet.   The choice done of the location 
of the interface points ensure that collisions occur one at a time so that the above  definition 
is unambiguous.   In general  there might be    triangles inside triangles.   
  The endpoints of the triangles 
  are    suitable coupled    pairs of interfaces points.    The  graphical representation  
just described   maps  each spin configuration in $ \XX_+ $ to a set of triangles.

  \vskip0.5cm
  \noindent {\bf Notation }
 {\it  Triangles will be usually denoted by $T$, the collection of triangles constructed as above 
  by $ \TT$ and we will write 
  $$ |T| = \hbox {cardinality of} \,\, \cap \Z =  \hbox {mass of} \quad T, \Eq (mars1)$$
 and by   $\hbox {supp} (T) \subset  \R $ the   basis of the triangle.}
 
 \medskip \noindent
  We have thus represented a configuration $\s \in \XX_+ $ as a collection of 
$\underline T = (T_1, \dots, T_n)$.  The above construction defines a one to
 one map from $ \XX_+ $ onto $ \TT$.  It is easy to see that a triangle 
configuration $\underline T$ belongs to $ \TT$ iff for any pair $T$ and 
$T'$ in $\underline T$ 
 $$ {\rm dist} (T,T')\ge \min\{ |T|, |T'|\} .\Eq (Ma1)$$ 
We say that  two collections of triangles $\underline S'$  and $\underline S$   are
 compatible  and we denote it by  $\underline S' \sim \underline S$  iff 
$ \underline S' \cup \underline S \in   \TT $ 
({\it i.e.} 
there exists a configuration in $\XX_+$ such that its corresponding collection  of triangles 
is the collection made of all triangles that are obtained by concatenating 
 $\underline S'$ and $\underline S$.)
By an abuse of notation, we write 
   $$  H^+_0(\underline T)= H^+_0(\sigma),  \quad G(\sigma(\underline T) )[\o]= G(\sigma)[\o], 
 \quad \s \in \XX_+ \iff \underline T \in   \TT $$. 
     \vskip0.5cm

  \noindent {\bf \Definition  (3)}   { \bf The energy difference } {\it 
 Given  two compatible   collections   of triangles  $\underline S \sim  \underline T$,    we denote
     $$
H^+(\underline S| \underline T):=   H^+(\underline S\cup \underline T)-H^+( \underline T). \Eq (D3) $$}

Let  $ \underline T= (T_1, \dots, T_n)$ with $|T_i| \le |T_{i+1}|$ then using \eqv(D3) one has
$$
H^+(\underline T)=H^+(T_1|\underline T\setminus T_1)+ H^+(\underline T\setminus T_1).
\Eq(D3bis)
$$

The following Lemma proved  in  [\rcite {CFMP}], see  Lemma 2.1 there,  gives a lower bound on the cost to ``erase'' 
triangles sequentially starting from the smallest ones.

\medskip \noindent 
 \noindent  { \bf \Lemma (CFMP1)   [\rcite {CFMP}] } {\it For $ \alpha \in [0, \frac {\ln3} {\ln2} -1)$   
and  $$\zeta= \zeta(\a)=  1 -2(2^\a-1)\Eq (eo1)$$ 
one has 
$$ H^+_0(T_1| \underline T \setminus  T_1 )  \ge \zeta |T_1|^\a.  \Eq (Ma2) $$ 
By iteration, 
for  
 any $1\le i \le n$ 
   $$ H^+_0(\cup_{\ell=1}^i T_{\ell} | \underline T \setminus [\cup_{\ell=1}^i T_{\ell} ]) 
 \ge \zeta \sum_{\ell=1}^i |T_\ell|^\a. \Eq (Ma2a) $$ 
 For $ \a=0$, \eqv (Ma2)  and \eqv (Ma2a) hold with $ |T_\ell|^\a$ replaced by $\log|T_\ell| + 4$.}
 \vskip0.5cm 
The  estimate \eqv (Ma2a)  involves contributions coming from the full set of triangles 
associated to a given spin configuration, starting from the triangle having the smallest 
mass.  To  implement a Peierls bound in our set up we need to ``localize'' the estimates 
 to compute the weight  of a  triangle or of a finite set of triangles in a generic configuration.  
In order to   do this  [\rcite {CFMP}]       introduced  the notion of contours as clusters of
 nearby triangles sufficiently far away from all other triangles.
  \medskip \noindent {\bf Contours} 
   A contour $ \G$ is a collection $\underline T$ of  triangles related by a 
hierarchical network of connections  controlled by a positive number $C$, see \eqv (SS1), under which all the triangles of a contour become mutually 
connected.     We denote by $T(\G)$ the triangle  whose basis is the smallest interval which contains all the 
triangles of the contour.  The right and left endpoints of $T(\G) \cap \Z$ are denoted by $x_{\pm} (\G)$.   
 We denote $|\G|$ the  mass of the contour  $\G$
 $$ |\G|= \sum_{T \in \G}|T| $$
 i.e.  $ |\G|$ is   the sum of the masses of all the triangles  belonging to $\G$.  
We denote  by   $ \RR (\cdot)$  the algorithm 
which 
associates to any configuration $ \underline T $  a configuration $ \{ \G_j\}$ of contours with the following 
properties.   
   \vskip0.5cm   \noindent
   { \bf P.0}  {\it   Let $ \RR (\underline T) = ( \G_1, \dots, \G_n)$, $ \G_i= \{ T_{j,i}, 1 \le j\le k_i\}$, 
then $\underline T=  \{ T_{j,i}, 1\le i \le n, 1 \le j\le k_i\}$} 
    
   \vskip0.5cm  \noindent
   { \bf P.1}  {\it Contours are well separated from each other.}  Any pair $ \G \neq \G'$ 
 verifies one of the following alternatives.
  
  $$       T (\G) \cap T (\G') = \emptyset $$
  i.e.  $[x_{-} (\G), x_{+} (\G)] \cap  [x_{-} (\G'), x_{+} (\G')] = \emptyset$,  in which case 
  $$ dist (\G, \G'):= \min_{T \in \G, T' \in \G'} dist ( T,T') > C  \min \left \{ |\G|^3, |\G'|^3\right \} \Eq (SS1)  $$
  where $C$ is a positive number.
 If 
   $$     T (\G) \cap T (\G')\neq  \emptyset, $$
   then either $T(\G)\subset T(\G')$ or $T(\G')\subset T(\G)$; moreover, supposing for instance
 that the former case is verified, (in which case we call $\G$ an inner contour) then for any 
triangle $ T'_i \in \G'$, either $T(\G)\subset T'_i$ or  $T(\G)\cap  T'_i= \emptyset $  and 
   $$ dist (\G, \G') >C |\G|^3, \quad \hbox {if} \quad  T(\G)\subset T(\G').  \Eq (SS3)$$
   \vskip0.5cm
   \noindent  
   { \bf P.2}  {\it Independence.}  Let $\{ \underline T^{(1)}, \dots, \underline T^{(k)}\}$,
 be $k>1$ configurations of triangles;    $\RR ( \underline T^{(i)}) = \{ \G_j^{(i)}, j=1,\dots, n_i\}$
 the contours of the configurations $\underline T^{(i)}$. Then if any distinct $ \G_j^{(i)}$ and 
 $\G_{j'}^{(i')}$ satisfies {\bf P.1}, 
   $$\RR ( \underline T^{(1)}, \dots, \underline T^{(k)} )= \{  \G_j^{(i)}, j=1,\dots, n_i; i=1, \dots, k \}. $$
  As proven in   [\rcite {CFMP}], the  algorithm $ \RR (\cdot)$ having properties  {\bf P.0},  {\bf P.1} and   {\bf P.2}   is unique  and therefore there is a 
bijection between  families of triangles  and contours. 
   Next we report the   estimates proven in  Theorem 3.2 of   [\rcite {CFMP}] which are essential  for  this paper. 
 \vskip0.5cm
  \noindent  { \bf \Theorem (1CFMP) [\rcite {CFMP}] } {\it  Let 
 $\a \in [0, \frac {\ln3} {\ln2} -1)$ and the constant $C$ given in  \eqv (SS1),    be so large that 
  $$ \sum_{m\ge 1} \frac {4m} {[C m]^3} \le \frac 12,\Eq (SS2)$$
  where  $ [x ]$ denotes  the integer part of $x$. 
For any $ \underline  T \in  \{ \underline T\}$,  let $\G_0\in \RR (\underline T)$ 
be a contour, $\underline S^{(0)}$ the triangles in $\G_0$ and     
$\zeta= \zeta(\a)=  1 -2(2^\a-1).$
Then    
$$H^+_0(  \underline S^{(0)}   | \underline T \setminus 
 \underline   S^{(0)}  )\ge \frac \zeta 2  \sum_{ T \in \underline  S^{(0)}    }  |T|^\a.   \Eq (fo.1) $$
 For $ \a =0$,  \eqv (fo.1)  holds with  $|T|^\a$  replaced by  $ \log |T| +4$ . }
\vskip0.5cm
 Next we summarize   the results of Theorem 4.1 of  [\rcite {CFMP}] stated for $ \a>0$ and the corresponding estimate for $\a=0$ given in Appendix F of   [\rcite {CFMP}]. 
 \vskip0.5cm

  \noindent  { \bf \Theorem  (2CFMP)  [\rcite {CFMP}]} {\it  For  any    $\a > 0$ 
  there exists  $C_0(\a)$ so that for $b\ge  C_0(\a)$  and for all $m>0$    
  $$ \sum_{\{0 \in \Gamma,  |\G| =m \}}   w_{b}^\a  (\Gamma)  \le 2m e^{- b m^{\a}},\Eq (E.6a) $$ 
where 
$$ w_b^\a (\Gamma):= \prod_{ T \in \G} e^{-b| T|^{\a}}.  \Eq (E.4a) $$ 
When $ \a  =0 $   
$$ w_b^0 (\Gamma):=  \prod_{ T \in \G}  e^{- b(\log |T| +4)} = \prod_{ T \in \G} \left ( |T|^{-b} e^{-4b} \right )    \Eq (E.4d) $$ 
and   there exists $C_0 $ so that for $b \ge  C_0 $ 
 $$ \sum_{\{0 \in \Gamma,  |\G| =m \}}   w_{b}^0  (\Gamma)  \le 2m e^{- b(\log m +4)}. \Eq (E.6ab) $$ }

\centerline{\bf References}
\vskip.3truecm
\item{[\rtag{ACCN}]} M. Aizenman, J. Chayes, L. Chayes and C. Newman:
{Discontinuity of the magnetization in one--dimensio\-nal $1/|x-y|^2$ percolation, Ising and Potts models.}
{\it J. Stat. Phys.} {\bf 50} no. 1-2 1--40 (1988).
\item{[\rtag{AW}]} M. Aizenman, and  J. Wehr:
{ Rounding of first or\-der pha\-se tran\-si\-tions
in sys\-tems with quenched disorder.}
{\it  Com. Math. Phys.} {\bf 130}, 489--528 (1990).
  \item{[\rtag{Bo}]} A. Bovier 
{\it  Statistical Mechanics of Disordered Systems.}
{ Cambridge Series in Statistical and Probabilistic mathematics.}, (2006).
  \item{[\rtag{BK}]} J. Bricmont, and A. Kupiainen.
{  Phase transition in the three-dimensional
random field Ising model.}
{\it Com. Math. Phys.},{\bf 116}, 539--572 (1988).
 \item{[\rtag {CFMP}]} M. Cassandro, P.  A. Ferrari, I. Merola and E. Presutti.  
{ Geometry of contours and Peierls estimates in $d=1$ Ising models with long range interaction.} 
{\it  J. Math. Phys.} {\bf 46}, no 5,    (2005) 
\item{[\rtag{COP1}]} M. Cassandro, E. Orlandi, and P.Picco.    
{ Typical configurations for one-dimensional random field Kac model. }
{\it   Ann. Prob.}  {\bf 27}, No 3, 1414-1467, (1999).
\item{[\rtag{COPV}]} M. Cassandro, E. Orlandi, P. Picco and
M.E. Vares. {One-dimensional random field Kac's model: Localization of the Phases}
{\it   Electron. J. Probab.}, {\bf 10}, 786-864,  (2005).  

 \item{[\rtag{COP3}]} M. Cassandro, E. Orlandi, and P.Picco.   
 { Phase Transition in the 1d Random Field   Ising Model  with  long range interaction. }
{\it   Comm. Math. Phy.}, {\bf 2},  731-744 (2009) 

  \item{[\rtag{CT}]}
Y.S. Chow and H. Teicher {\it Probability theory. Independence, interchangeability, martingales.}
 Third edition. Springer Texts in Statistics. Springer-Verlag, New York, 1997. 

\item{[\rtag{D0}]} R. Dobrushin: 
The description of a random field by means of conditional probabilities and. conditions of its regularity.
{\it  Theory Probability Appl.} {\bf 13}, 197-224 (1968)  
\item{[\rtag{D1}]} R. Dobrushin:
The conditions of absence of phase transitions in one-dimensional classical systems.
{\it Matem. Sbornik}, {\bf 93} (1974), N1, 29-49

\item{[\rtag{D2}]} R. Dobrushin:
Analyticity of correlation functions in one-dimensional classical systems with slowly decreasing potentials.
{\it Comm. Math. Phys.} {\bf 32} (1973), N4, 269-289

  \item{[\rtag{Dy}]} F.J. Dyson:
 {  Existence of  phase transition in a  one-dimensional Ising ferromagnetic.} 
{ \it Comm. Math. Phys.},{\bf 12},91--107,  (1969).

\item{[\rtag{FS}]} J.  Fr\"ohlich and  T. Spencer:  
 {  The phase transition in the one-dimensional Ising model with
 $\frac 1 {r^2}$ interaction energy.} { \it Comm. Math. Phys.}, {\bf 84}, 87--101,   (1982).

\item{[\rtag{I}]} J.Z. Imbrie 
{ Decay of correlations in the one-dimensional Ising model with $J_{ij}=\mid i-j\mid^{-2}$,} 
{\it Comm. Math. Phys.} {\bf  85}, 491--515. (1982).
 
\item{[\rtag{IN}]} J.Z. Imbrie and C.M. Newman 
{ An intermediate phase with slow decay of correlations in one-dimensional
 $1/\vert x-y\vert ^2$ percolation, Ising and Potts models.}
 {\it  Comm. Math. Phys.}  {\bf 118} , 303--336 (1988). 

\item{[\rtag{KAH}]} J. P. Kahane Propri\'et\'es locale des fonctions \`a s\'eries de 
Fourier al\'eatoires {\it Studia Matematica} {\bf 19} 1--25 (1960).

 \item{[\rtag{Lcam}]} L.  Le Cam   
 { \it Asymptotic methods in statistical decision theory.} { Springer Series in Statistics},
 Springer-Verlag, New York, Berlin, Heidelberg,  (1986).

 \item{[\rtag{OP}]}  E. Orlandi, and P.Picco.    { One-dimensional random field  Kac's model:
 weak large deviations principle. } 
   { \it Electronic Journal Probability}  {\bf 14},   1372--1416, (2009).

\item {[\rtag{RT}]} J. B.  Rogers and C.J. Thompson:  
Absence of long range order in one dimensional spin systems. 
{\it  J. Statist. Phys.} {\bf 25}, 669--678 (1981)

\item {[\rtag{Ru}]} D. Ruelle:  
Statistical mechanics of one-dimensional Lattice gas.
{\it  Comm. Math. Phys.} {\bf 9}, 267--278 (1968)

  \end